\documentclass[nosumlimits,twoside,reqno]{amsart}
\usepackage{amssymb,amsmath,amsfonts}
\usepackage[bookmarksnumbered,plainpages,hypertex]{hyperref}
\usepackage{float}
\usepackage{graphicx}

\newtheorem{theorem}{\sc Theorem}[section]
\newtheorem{proposition}[theorem]{\sc Proposition}
\newtheorem{lemma}[theorem]{\sc Lemma}
\newtheorem{corollary}[theorem]{\sc Corollary}
\theoremstyle{definition}
\newtheorem{definition}[theorem]{\sc Definition}

\newtheorem{example}[theorem]{\sc Example}

\newtheorem{noname}[theorem]{}
\theoremstyle{remark}
\newtheorem{remark}[theorem]{\sc Remark}
\newtheorem{remarks}[theorem]{\sc Remarks}

\newtheoremstyle{mytheorem}{}{}{\it}{0pt}{\bf}{}{0pt}
{\sc \thmname{#1}\ \thmnote{#3}\ } \theoremstyle{mytheorem}

\begin{document}
\title{A Milnor-Moore Type Theorem for Braided Bialgebras}
\author{A. Ardizzoni, C. Menini}
\author{D. \c{S}tefan}
\date{}

\begin{abstract}
The paper is devoted to prove a version of Milnor-Moore Theorem
for connected braided bialgebras that are infinitesimally
cocommutative. Namely in characteristic different from $2$, we
prove that, for a given connected braided bialgebra $A$ having a
$\lambda $-cocommutative infinitesimal braiding for some regular
element $\lambda \neq 0$ in the base field, then the infinitesimal
braiding of $A$ is of Hecke-type of mark $\lambda $ and $A$ is
isomorphic as a braided bialgebra to the symmetric algebra of the
braided subspace of its primitive elements.
\end{abstract}

\thanks{This paper was written while A. Ardizzoni and C. Menini were members
of GNSAGA with partial financial support from MiUR. D. \c{S}tefan
was partially supported by INDAM, while he was visiting professor
at University of Ferrara, and by CERES Project 4-147/2004.}
\maketitle
%\tableofcontents

\section*{Introduction}

The structure of cocommutative connected bialgebras is well-understood in
characteristic zero. By Milnor-Moore Theorem \cite{MM} such a bialgebra is
the enveloping algebra of its primitive part, regarded as a Lie algebra in a
canonical way. This result is one of the most important ingredients in the
proof of Cartier-Gabriel-Kostant Theorem \cite{Di, Sw}, that characterizes
cocommutative pointed bialgebras in characteristic zero as \textquotedblleft
products\textquotedblright\ between enveloping algebras and group algebras.

Versions of Milnor-Moore Theorem and Cartier-Gabriel-Kostant Theorem for
graded bialgebras (over $\mathbb{Z}$ and $\mathbb{Z}_{2}$) can be also found
in the work of Kostant, Leray and Milnor-Moore. For other more recent
results, analogous to Milnor-Moore Theorem, see \cite{Go, LR, R1, R2}.

The $\mathbb{Z}_{2}$-graded bialgebras, nowadays called superbialgebras,
appeared independently in the work of Milnor-Moore and MacLane. From a
modern point of view, superalgebras can be seen as bialgebras in a braided
monoidal category. These structures play an increasing role not only in
algebra (classification of finite dimensional Hopf algebras or theory of
quantum groups), but also in other fields of mathematics (algebraic
topology, algebraic groups, Lie algebras, etc.) and physics. Braided
monoidal categories were formally defined by Joyal and Street in the seminal
paper \cite{JS}, while (bi)algebras in a braided category were introduced in
\cite{Ma}. By definition, $(A,\nabla ,u,\Delta ,\varepsilon )$ is a
bialgebra in a braided category $\mathfrak{M,}$ if $\nabla $ is an
associative multiplication on $A$ with unit $u$ and $\Delta $ is a
coassociative comultiplication on $A$ with counit $\varepsilon $ such that $%
\varepsilon $ and $\Delta $ are morphism of algebras where the
multiplication on $A\otimes A$ is defined via the braiding. In
other words, the last property of $\Delta $ rereads as follows
\begin{equation}
\Delta \nabla =(\nabla \otimes \nabla )(A\otimes c_{A,A}\otimes A)(\Delta
\otimes \Delta ),  \label{introd}
\end{equation}%
where $c_{X,Y}:X\otimes Y\rightarrow Y\otimes X$ denotes the braiding in $%
\mathfrak{M}.$ The relation above had already appeared in a natural way in
\cite{Rad}, where Hopf algebras with a projection are characterized. More
precisely, let $A$ be a Hopf algebra and let $p:A\rightarrow A$ be a
morphism of Hopf algebras such that $p^{2}=p$. To these data, Radford
associates an ordinary Hopf algebra $B:=\mathrm{Im}(p)$ and a Hopf algebra $%
R=\{a\in A\mid (A\otimes p)\Delta (a)=a\otimes 1\}$ in the braided category $%
_{B}^{B}\mathcal{YD}$ of Yetter-Drinfeld modules. Then he shows that $%
A\simeq R\otimes B$, where on $R\otimes B$ one puts the tensor product
algebra and tensor product coalgebra of $R$ and $B$ in the category $_{B}^{B}%
\mathcal{YD}$ (the braiding of $_{B}^{B}\mathcal{YD}$ is used to twist the
elements of $R$ and $B$). It is worth to notice that the above Hopf algebra
structure on $R\otimes B$ can be constructed for an arbitrary Hopf algebra $R
$ in $_{B}^{B}\mathcal{YD}.$ It is called the \emph{bosonization} of $R$ and
it is denoted by $R\#B$. Notably, the \textquotedblleft
product\textquotedblright\ that appears in Cartier-Gabriel-Kostant Theorem
is precisely the bosonization of an enveloping Lie algebra, regarded as a
bialgebra in the category of Yetter-Drinfeld modules over its coradical. The
result of Radford was extended for more general classes of bialgebras with a
projection in \cite{AMS2, AMS3}.

Bosonization also plays a very important role in the \textquotedblleft
lifting\textquotedblright\ method for the classification of finite
dimensional pointed Hopf algebras, introduced by N. Andruskiewitsch and H.J.
Schneider. Roughly speaking, the lifting method requires two steps. If $A$
is a pointed Hopf algebra, then $\mathrm{gr}A$, the graded associated of $A$
with respect to the coradical filtration, is a Hopf algebra with projection
onto the homogeneous component of degree $0.$ Hence, by Radford's result, $%
\mathrm{gr}A$ is the bosonization of a graded \emph{connected} Hopf algebra $%
R$ in ${}_{B}^{B}\mathcal{YD}$, where $B$ is the coradical of $A$.
Accordingly to the lifting method, first one has to classify all
connected and graded Hopf algebras $R$ in $_{B}^{B}\mathcal{YD}$
such that $\dim
R=\dim A/\dim B$. Then one has to find all Hopf algebras $A$ such that $%
A\simeq R\#B$, with $R$ as in the first step.

Therefore, in many cases, for proving a certain property of Hopf algebras,
it is enough to do it in the connected case. The price that one has to pay
is that we have to work with Hopf algebras in a braided category (usually$%
{}_{B}^{B}\mathcal{YD}$), and not with ordinary Hopf algebras. Motivated by
this observation, in this paper we will investigate connected and
cocommutative bialgebras in a braided category. Actually nowadays people
recognize that it is more appropriate to work with braided bialgebras, that
were introduced in \cite{Ta} (see e.g. \cite{Gu- the YB} and \cite%
{Kharchenko- connected}).

To define a braided bialgebra we first need a braided vector space, that is
a pair $(A,\mathfrak{c})$, where $A$ is a vector space and $\mathfrak{c}$ is
a solution of the braid equation (\ref{ec: braid equation}). Then we need an
algebra $(A,\nabla ,u)$ and a coalgebra $(A,\Delta ,\epsilon )$ which are
compatible with the braiding (see Definitions \ref{def: braided algebra} and %
\ref{def: braided coalgebra}). Now, for defining braided bialgebras, one can
proceed as in the classical case; see \ref{de: c-bialgebra}.

The prototype braided bialgebra is $T:=T(V,\mathfrak{c}),$ the tensor
algebra of a braided $K$-vector space $(V,\mathfrak{c}).$ The braiding $%
\mathfrak{c}$ lifts uniquely to an operator $\mathfrak{c}_{T}:T\otimes
T\rightarrow T\otimes T$. The usual algebra structure on $T(V)$ is
compatible with $\mathfrak{c}_{T}$, so $T\otimes T$ is an algebra (for the
multiplication, of course, we use $\mathfrak{c}_{T}$ and not the usual flip
map). Therefore, there is a unique coalgebra structure on $T$, so that $V$
is included in the space of primitive elements of $T$ and the
comultiplication is an algebra map.

For constructing other examples of braided bialgebras, we focus on the case
when $(V,\mathfrak{c})$ is a braided vector space such that $\mathfrak{c}$
is a braiding of Hecke-type of mark $\lambda \in K^{\ast }$, that is $%
\mathfrak{c}$ is a root of $(X+1)(X-\lambda ).$ Then, to every $K$-linear
map $\mathfrak{b}:V\otimes V\rightarrow V$ which is compatible with $%
\mathfrak{c}{}$(i.e. a so called braided bracket, see Definition \ref{ec: b
compatibilitate}) we associate a new braided bialgebra $U(V,\mathfrak{c},%
\mathfrak{b}),$ called enveloping algebra, as follows. The set
\begin{equation*}
X_{\mathfrak{c},\mathfrak{b}}=\{\mathfrak{c}(z)-\lambda z-\mathfrak{b}%
(z)\mid z\in V\otimes V=T^{2}(V)\}
\end{equation*}%
contains only primitive elements in $T(V,\mathfrak{c})$, so the ideal $I_{%
\mathfrak{c},\mathfrak{b}}$ generated by $X_{\mathfrak{c},\mathfrak{b}}$ in $%
T(V,\mathfrak{c})$ is a coideal too. Hence the quotient $U(V,\mathfrak{c},%
\mathfrak{b})$ of $T(V,\mathfrak{c})$ through $I_{\mathfrak{c},\mathfrak{b}}$
is a braided bialgebra. As a particular case we obtain the braided symmetric
algebra $S(V,\mathfrak{c}):=U(V,\mathfrak{c},0)$.

It is worthwhile noticing that the braided subspace $(P,\mathfrak{c}_{P})$
of primitive elements of a connected braided bialgebra $(A,\mathfrak{c}_{A})$
can always be endowed with a braided bracket $\mathfrak{b}_{P}:P\otimes
P\rightarrow P$ whenever $\mathfrak{c}_{P}$ is a braiding of Hecke-type. In
this case one can consider the corresponding enveloping algebra (see Theorem
(\ref{pr: universalitate U})).

The main result of the paper is Theorem \ref{te: MM}. It is a version of
Milnor-Moore Theorem for connected braided bialgebras that are
infinitesimally cocommutative. In fact, we prove that, for a given field $K$
with $\mathrm{char}\,K\neq 2$, if $A$ is a connected braided bialgebra
having a $\lambda $-cocommutative infinitesimal braiding $\mathfrak{c}%
_{P\left( A\right) }$ for some regular element $\lambda \neq 0$ in $K$, then

\begin{itemize}
\item the infinitesimal braiding of $A$ is of Hecke-type of mark $\lambda $
and

\item $A$ is isomorphic as a braided bialgebra to the symmetric algebra $%
S\left( P\left( A\right) ,\mathfrak{c}_{P\left( A\right) }\right) $ of $%
\left( P\left( A\right) ,\mathfrak{c}_{P\left( A\right) }\right) $ whenever $%
\lambda \neq 1$.
\end{itemize}

We point out that the case $\lambda =1$ can be easily treated by means of
Kharchenko's results in \cite{Kharchenko- connected} (see Remark \ref{rem:
MM}).

To achieve our result we characterize bialgebras of type one with
infinitesimal braiding of Hecke type (see Theorem \ref{teo: aureo}).
Moreover, for a given braided vector space $(V,\mathfrak{c})$ of Hecke-type
of mark $\lambda \neq 0,1$, we show that $\mathfrak{b}=0$ is the unique $%
\mathfrak{c}$-bracket on $(V,\mathfrak{c})$ for which the $K$-linear
canonical map $\iota _{\mathfrak{c},\mathfrak{b}}:V\rightarrow U(V,\mathfrak{%
c},\mathfrak{b})$ is injective whenever $(3)!_{\lambda }\neq 0$ (see Theorem %
\ref{te: Lie algebras}).

\section{Braided bialgebras}

Throughout this paper $K$ will denote a field. All vector spaces will be
defined over $K$ and the tensor product of two vector spaces will be denoted
by $\otimes$.

In this section we define the main notion that we will deal with, namely
braided bialgebras. We also introduce one of the basic examples, namely the
tensor algebra of a braided vector space.

\begin{definition}
A pair $(V,\mathfrak{c})$ is called \emph{braided vector space} if $%
\mathfrak{c}:V\otimes V\rightarrow V\otimes V$ is a solution of the braid
equation
\begin{equation}
\mathfrak{c}_{1}\mathfrak{c}_{2}\mathfrak{c}_{1}=\mathfrak{c}_{2}\mathfrak{c}%
_{1}\mathfrak{c}_{2}  \label{ec: braid equation}
\end{equation}%
where $\mathfrak{c}_{1}=\mathfrak{c}\otimes V$ and $\mathfrak{c}%
_{2}=V\otimes \mathfrak{c}$. A morphism of braided vector spaces $(V,%
\mathfrak{c}_{V})$ and $(W,\mathfrak{c}_{W})$ is a $K$-linear map $%
f:V\rightarrow W$ such that $\mathfrak{c}_{W}(f\otimes f)=(f\otimes f)%
\mathfrak{c}_{V}.$
\end{definition}

Note that, for every braided vector space $(V,\mathfrak{c})$ and every $%
\lambda \in K,$ the pair $(V,\lambda \mathfrak{c})$ is a braided vector
space too. A general method for producing braided vector spaces is to take
an arbitrary braided category $(\mathcal{M{}},\otimes ,K,a,l,r,c),$ which is
a monoidal subcategory of the category of $K{}$-vector spaces. Hence any
object $V\in \mathcal{M}$ can be regarded as a braided vector space with
respect to $\mathfrak{c}:=c_{V,V}.$ Here, $c_{X,Y}:X\otimes Y\rightarrow
Y\otimes X$ denotes the braiding in $\mathcal{M{}}.$ The category of
comodules over a coquasitriangular Hopf algebra and the category of
Yetter-Drinfeld modules are examples of such categories. More particularly,
every bicharacter of a group $G$ induces a braiding on the category of $G$
-graded vector spaces.

\begin{definition}[Baez, \protect\cite{Ba}]
\label{def: braided algebra}A \emph{braided algebra}, or $\mathfrak{c}$-%
\emph{algebra} is a quadruple $(A,\nabla ,1,\mathfrak{c})$ where $(A,%
\mathfrak{c})$ is a braided vector space and $(A,\nabla ,1)$ is an
associative unital algebra such that $\nabla $ and $u$ commute with $%
\mathfrak{c}$, that is the following conditions hold:
\begin{gather}
\mathfrak{c}(\nabla \otimes A)=(A\otimes \nabla )(\mathfrak{c}\otimes
A)(A\otimes \mathfrak{c}),  \label{ec: c-alg1} \\
\mathfrak{c}(A\otimes \nabla )=(\nabla \otimes A)(A\otimes \mathfrak{c})(%
\mathfrak{c}\otimes A),  \label{ec: c-alg2} \\
\mathfrak{c}(1\otimes a)=a\otimes 1,\qquad \mathfrak{c}(a\otimes 1)=1\otimes
a,\qquad \forall a\in A.  \label{ec: c-alg3}
\end{gather}%
A morphism of braided algebras is, by definition, a morphism of ordinary
algebras which, in addition, is a morphism of braided vector spaces.
\end{definition}

\begin{remark}
\label{re: nabla} 1) Let $(A,\nabla ,u,\mathfrak{c})$ is a braided algebra.
Then $A\otimes A$ is an associative algebra with multiplication $\nabla
_{A\otimes A}:=(\nabla \otimes \nabla )(A\otimes \mathfrak{c}\otimes A)$ and
unit $1\otimes 1.$ Moreover, $A\otimes A$ is a $\mathfrak{c}_{A\otimes A}$%
-algebra, where $\mathfrak{c}_{A\otimes A}=(A\otimes \mathfrak{c}\otimes A)(%
\mathfrak{c}\otimes \mathfrak{c})(A\otimes \mathfrak{c}\otimes A)$. This
algebra structure on $A\otimes A$ will be denoted by $A\otimes _{\mathfrak{c}%
}A$.

2)If $A$ is an object in a braided monoidal category $\mathcal{M{}}$ and $%
\mathfrak{c}:=c_{A,A}$ then the above four compatibility relations hold
automatically, as the braiding $c$ is a natural morphism.
\end{remark}

\begin{definition}
\label{def: braided coalgebra}A \emph{braided coalgebra} (or $\mathfrak{c}$-%
\emph{coalgebra}) is a quadruple $(C,\Delta ,\varepsilon ,\mathfrak{c})$
where $(C,\mathfrak{c})$ is a braided vector space and $(C,\Delta
,\varepsilon )$ is a coassociative counital coalgebra such that the
comultiplication $\Delta $ and the counit $\varepsilon $ commute with $%
\mathfrak{c}$, that is the following relations hold:
\begin{gather}
(C\otimes \Delta )\mathfrak{c}=(\mathfrak{c}\otimes C)(C\otimes \mathfrak{c}%
)(\Delta \otimes C),  \label{c-c1} \\
(\Delta \otimes C)\mathfrak{c}=(C\otimes \mathfrak{c})(\mathfrak{c}\otimes
C)(C\otimes \Delta ),  \label{c-c2} \\
(\varepsilon \otimes C)\mathfrak{c}(c\otimes d)=\varepsilon (d)c=(C\otimes
\varepsilon )\mathfrak{c}(d\otimes c),\qquad \forall c,d\in C.  \label{c-c3}
\end{gather}%
A morphism of braided coalgebras is, by definition, a morphism of ordinary
coalgebras which, in addition, is a morphism of braided vector spaces.
\end{definition}

\begin{noname}
\label{nn: connected} Recall that a coalgebra $C$ is called \emph{connected}
if the coradical $C_{0}$ of $C$ is one dimensional. In this case there is a
unique group-like element $g\in C$ such that $C_{0}=Kg$. Sometimes, we will
write $(C,g)$, to emphasize the group-like element $g$. We also ask that $%
f(g_{C})=g_{D}$, for any morphism $f:(C,g_{C})\rightarrow (D,g_{D})$ of
connected coalgebras.

By definition, a $\mathfrak{c}$-coalgebra $C$ is \emph{connected} if $%
C_{0}=Kg$ and, for any $x\in C$,
\begin{equation}  \label{de: connected}
\mathfrak{c}(x\otimes g)=g\otimes x,\qquad \mathfrak{c} (g\otimes
x)=x\otimes g.
\end{equation}
\end{noname}

\begin{definition}[Takeuchi, \protect\cite{Ta}]
\label{de: c-bialgebra}A \emph{braided bialgebra} is a sextuple $(A,\nabla
,1,\Delta ,\varepsilon ,\mathfrak{c})$ where

\begin{itemize}
\item $(A,\nabla ,1,\mathfrak{c})$ is a braided algebra,

\item $(A,\Delta ,\varepsilon ,\mathfrak{c})$ is a braided coalgebra,

\item $\Delta $ and $\varepsilon $ are morphisms of algebras (on the vector
space $A\otimes A$ we take the algebra structure $A\otimes _{\mathfrak{c}}A$%
).
\end{itemize}
\end{definition}

\begin{remark}
Note that $\Delta :A\rightarrow A\otimes _{\mathfrak{c}}A$ is multiplicative
if and only if
\begin{equation}
\Delta \nabla =(\nabla \otimes \nabla )(A\otimes \mathfrak{c}{}\otimes
A)(\Delta \otimes \Delta ).  \label{ec: c-bialgebra}
\end{equation}
\end{remark}

\begin{noname}
We will need graded versions of braided algebras, coalgebras and bialgebras.
By definition, a braided algebra $(A,\nabla ,1,\mathfrak{c})$ is graded if $%
A=\bigoplus_{n\in \mathbb{N}}A^{n}$ and $\nabla (A^{n}\otimes
A^{m})\subseteq A^{n+m}$. The braiding $\mathfrak{c}$ is assumed to satisfy $%
\mathfrak{c}(A^{n}\otimes A^{m})\subseteq A^{m}\otimes A^{n}.$ In this case,
it is easy to see that $1\in A^{0}$.

Therefore a graded braided algebra can be defined by means of maps $\nabla
^{n,m}:A^{n}\otimes A^{m}\rightarrow A^{n+m}$ and $\mathfrak{c}%
^{n,m}:A^{n}\otimes A^{m}\rightarrow A^{m}\otimes A^{n}$, and an element $%
1\in A^{0}$ such that:
\begin{gather}
\nabla ^{n+m,p}(\nabla ^{n,m}\otimes A^{p})=\nabla ^{n,m+p}(A^{n}\otimes
\nabla ^{m,p}),\quad n,m,p\in \mathbb{N},  \label{gr1} \\
\nabla ^{0,n}(1\otimes a)=a=\nabla ^{n,0}(a\otimes 1),\quad \forall a\in
A^{n},\quad n\in \mathbb{N}.  \label{gr2} \\
\mathfrak{c}^{n+m,p}(\nabla ^{n,m}\otimes A^{p})=(A^{p}\otimes \nabla
^{n,m})(\mathfrak{c}^{n,p}\otimes A^{m})(A^{n}\otimes \mathfrak{c}^{m,p}{}),
\label{gr3} \\
\mathfrak{c}^{n,m+p}(A^{n}\otimes \nabla ^{m,p})=(\nabla ^{m,p}\otimes
A^{n})(A^{m}\otimes \mathfrak{c}^{n,p})(\mathfrak{c}^{n,m}\otimes A^{p}),
\label{gr4} \\
\mathfrak{c}^{0,n}(1\otimes a)=a\otimes 1\qquad \text{and\qquad }\mathfrak{c}%
^{n,0}(a\otimes 1)=1\otimes a,\qquad \forall a\in A^{n}.  \label{gr5}
\end{gather}%
The multiplication $\nabla $ can be recovered from $(\nabla ^{n,m})_{n,m\in
\mathbb{N}}$ as the unique $K$-linear map such that $\nabla (x\otimes
y)=\nabla ^{p,q}(x\otimes y),\forall p,q\in \mathbb{N},\forall x\in
A^{p},\forall y\in A^{q}.$ Analogously, the braiding $\mathfrak{c}$ is
uniquely defined by $\mathfrak{c}(x\otimes y)=\mathfrak{c}^{p,q}(x\otimes
y),\forall p,q\in \mathbb{N},\forall x\in A^{p},\forall y\in A^{q}.$ We will
say that $\nabla ^{n,m}$ and $\mathfrak{c}^{n,m}$ are the $(n,m)$%
-homogeneous components of $\nabla $ and $\mathfrak{c}$, respectively.

Graded braided coalgebras can by described in a similar way. By definition a
braided coalgebra $(C,\Delta ,\varepsilon ,\mathfrak{c}{})$ is graded if $%
C=\bigoplus_{n\in \mathbb{N}}C^{n},\Delta (C^{n})\subseteq
\sum_{p+q=n}C^{p}\otimes C^{q}$, $\mathfrak{c}(C^{n}\otimes C^{m})\subseteq
C^{m}\otimes C^{n}$ and $\varepsilon |_{C_{n}}=0$, for $n>0$ . If $\pi ^{p}$
denotes the projection onto $C^{p}$ then the comultiplication $\Delta $ is
uniquely defined by the maps $\Delta ^{p,q}:C^{p+q}\rightarrow C^{p}\otimes
C^{q}$, where $\Delta ^{p,q}:=(\pi ^{p}\otimes \pi ^{q})\Delta |_{C^{p+q}}$.
The counit is given by a map $\varepsilon ^{0}:C^{0}\rightarrow K,$ while
the braiding $\mathfrak{c}{}$ is uniquely determined by a family $(\mathfrak{%
c}^{n,m})_{n,m\in \mathbb{N}},$ as for braided algebras. The families $%
(\Delta ^{n,m})_{n,m\in \mathbb{N}},$ $(\mathfrak{c}^{n,m})_{n,m\in \mathbb{N%
}}$ and $\varepsilon ^{0}$ has to satisfy the relations that are dual to (%
\ref{gr1}) -- (\ref{gr5}), namely:
\begin{gather}
(\Delta ^{n,m}\otimes C^{p})\Delta ^{n+m,p}=(C^{n}\otimes \Delta
^{m,p})\Delta ^{n,m+p},\quad n,m,p\in \mathbb{N},  \label{c1} \\
(\varepsilon ^{0}\otimes C^{n})\Delta ^{0,n}(c)=c=(C^{n}\otimes \varepsilon
^{0})\Delta ^{n,0}(c),\quad \forall c\in C^{n},\quad n\in \mathbb{N}.
\label{c2} \\
{}(C^{p}\otimes \Delta ^{n,m})\mathfrak{c}^{n+m,p}=(\mathfrak{c}%
^{n,p}\otimes C^{m}{})(C^{n}\otimes \mathfrak{c}^{m,p})(\Delta ^{n,m}\otimes
C^{p}),  \label{c3} \\
(\Delta ^{m,p}\otimes C^{n})\mathfrak{c}^{n,m+p}=(C^{m}\otimes \mathfrak{c}%
^{n,p})(\mathfrak{c}^{n,m}\otimes C^{p})(C^{n}\otimes \Delta ^{m,p}),
\label{c4} \\
(\varepsilon ^{0}\otimes C)\mathfrak{c}(c\otimes d)=\varepsilon
^{0}(d)c=(C\otimes \varepsilon ^{0})\mathfrak{c}(d\otimes c),\qquad \forall
c\in C^{n},\ \ \forall d\in C^{0}.  \label{c5}
\end{gather}%
We will say that $\Delta ^{n,m}$ is the $(n,m)$-homogeneous component of $%
\Delta $.

A \textbf{graded braided bialgebra} is a braided bialgebra which is graded
both as an algebra and as a coalgebra.
\end{noname}

\begin{remark}
\label{re: graded connected} Let $C=\bigoplus_{n\in \mathbb{N}}C^{n}$ be a
graded braided coalgebra. By \cite[Proposition 11.1.1]{Sw}, if $%
(C_{n})_{n\in \mathbb{N}}$ is the coradical filtration, then $C_{n}\subseteq
\bigoplus_{m\leq n}C^{m}$. Therefore, if $C^{0}$ is one dimensional then $C$
is connected.
\end{remark}

\begin{definition}
A graded braided coalgebra will be called $0$\emph{-connected} if its
homogeneous component of degree $0$ is of dimension one.
\end{definition}

\begin{lemma}
\label{le: gr(C)} Let $(C,\mathfrak{c})$ be a connected braided coalgebra.
Then $\mathfrak{c}$ induces a canonical braiding $\mathfrak{c}_{\mathrm{gr\,}%
C}$ on $\mathrm{gr\,}C$ such that $(\mathrm{gr\,}C,\mathfrak{c}_{\mathrm{gr\,%
}C})$ is a $0$-connected graded braided coalgebra, where $\mathrm{gr\,}C$ is
constructed with respect to the coradical filtration on $C.$
\end{lemma}

\begin{proof}
Let $(C_{n})_{n\in \mathbb{N}}$ be the coradical filtration. Since $C$ is
connected, we have $C_{0}=Kg$, where $g$ is the unique group-like element of
$C$. We claim that $\mathfrak{c}(C_{n}\otimes C_{m})\subseteq C_{m}\otimes
C_{n}.$ For $n=0$ this relation holds true as, by definition, $\mathfrak{c}%
(x\otimes g)=g\otimes x,$ $\forall x\in C.$ We choose a basis $\{y_{i}\mid
i\in I\}$ on $C_{m}$ and we assume that the above inclusion is true for $n.$
Let $x\in C_{n+1}$ and $y\in C_{m}.$ By Taft-Wilson Theorem $\Delta
(x)=x\otimes g+g\otimes x+\sum_{k=1}^{q}x_{k}^{\prime }\otimes x_{k}^{\prime
\prime },$ where $x_{k}^{\prime },x_{k}^{\prime \prime }\in C_{n}.$ Moreover,%
\begin{equation*}
\mathfrak{c}(x\otimes y)=\sum_{i\in I}y_{i}\otimes x_{i},
\end{equation*}%
with $x_{i}\in C,$ and the set $\{i\in I\mid x_{i}\neq 0\}$ finite. Since $%
(C\otimes \Delta )\mathfrak{c}=(\mathfrak{c}\otimes C)(C\otimes \mathfrak{c}%
)(\Delta \otimes C)$ we get
\begin{equation*}
\sum_{i=1}^{p}y_{i}\otimes \Delta (x_{i})=\sum_{i=1}^{p}y_{i}\otimes
x_{i}\otimes g+\sum_{i=1}^{p}y_{i}\otimes g\otimes x_{i}+(\mathfrak{c}%
\otimes C)(C\otimes \mathfrak{c})(\sum_{k=1}^{q}x_{k}^{\prime }\otimes
x_{k}^{\prime \prime }\otimes y).
\end{equation*}%
By induction hypothesis, $(\mathfrak{c}\otimes C)(C\otimes \mathfrak{c}%
)(\sum_{k=1}^{q}x_{k}^{\prime }\otimes x_{k}^{\prime \prime }\otimes y)\in
C_{m}\otimes C_{n}\otimes C_{n},$ so this element can be written as $%
\sum_{i\in I}y_{i}\otimes z_{i},$ with $z_{i}\in C_{n}\otimes C_{n}$ and the
set $\{i\in I\mid z_{i}\neq 0\}$ is finite. Hence, for all $i\in I,$ we have
\begin{equation*}
\Delta (x_{i})=x_{i}\otimes g+g\otimes x_{i}+z_{i}.
\end{equation*}%
Thus $x_{i}\in C_{n+1},$ so $\mathfrak{c}(x\otimes y)\in C_{m}\otimes
C_{n+1}.$ Hence, by induction, $\mathfrak{c}(C_{n}\otimes C_{m})$ $\subseteq
C_{m}\otimes C_{n},$ so%
\begin{equation*}
\mathfrak{c}(C_{n-1}\otimes C_{m}+C_{n}\otimes C_{m-1})\subseteq
C_{m}\otimes C_{n-1}+C_{m-1}\otimes C_{n}.
\end{equation*}%
Therefore $\mathfrak{c}$ induces a unique $K$-linear map $\mathfrak{c}_{%
\mathrm{gr\,}C}^{n,m}:\mathrm{gr\,}^{n}C\otimes \mathrm{gr\,}%
^{m}C\rightarrow \mathrm{gr\,}^{m}C\otimes \mathrm{gr\,}^{n}C.$ We define $%
\mathfrak{c}_{\mathrm{gr\,}C}:=\bigoplus_{n,m}\mathfrak{c}_{\mathrm{gr\,}%
C}^{n,m}.$ Now it is easy to see that $(\mathrm{gr\,}C,\mathfrak{c}_{\mathrm{%
gr\,}C})$ is a graded braided coalgebra.
\end{proof}

\begin{remark}
\label{rem: c_P}If $(C,\mathfrak{c},g)$ is a connected braided coalgebra
then $\mathfrak{c}$ induces a canonical braiding $\mathfrak{c}_{P}$ on the
space $P(C)=\{c\in C\mid \Delta (c)=c\otimes g+g\otimes c\}$, of primitives
elements in $C$. Indeed, by Taft-Wilson Theorem we have $P(C)=\left( \mathrm{%
Ker}\,\varepsilon \right) \cap C_{1}$. Thus $\mathfrak{c}$ maps $P(C)\otimes
P(C)$ to itself, see the proof of the preceding lemma.
\end{remark}

\begin{lemma}
\label{pr: Produs}Let $(A,\nabla ,u,\mathfrak{c})$ be a $\mathfrak{c}$%
-algebra and let $\Delta :A\rightarrow A\otimes _{\mathfrak{c}}A$ be a
morphism of algebras. We fix $x,a\in A$ such that $\Delta (a)=a\otimes
1+1\otimes a$. Then
\begin{equation}
\Delta (xa)=(\nabla \otimes A)(A\otimes \mathfrak{c})[\Delta (x)\otimes
a]+(A\otimes \nabla )[\Delta (x)\otimes a].  \label{ec: Produs1}
\end{equation}
\end{lemma}

\begin{proof}
Let $\Delta (x)=\sum_{i=1}^{p}x_{i}^{\prime }\otimes x_{i}^{\prime \prime }.$
Then%
\begin{eqnarray*}
\Delta (xa) &=&\Delta (x)\cdot \left( a\otimes 1+1\otimes a\right)
=\sum_{i=1}^{p}\left( x_{i}^{\prime }\otimes x_{i}^{\prime \prime }\right)
\cdot \left( a\otimes 1+1\otimes a\right) \\
&=&\sum_{i=1}^{p}x_{i}^{\prime }\mathfrak{c}(x_{i}^{\prime \prime }\otimes
a)+\sum_{i=1}^{p}x_{i}^{\prime }\otimes x_{i}^{\prime \prime }a \\
&=&(\nabla \otimes A)(A\otimes \mathfrak{c})[\Delta (x)\otimes a]+(A\otimes
\nabla )[\Delta (x)\otimes a].
\end{eqnarray*}
\end{proof}

\begin{proposition}
\label{pr: Delta(c)}To any braided vector space $(V,\mathfrak{c})$ we can
associate a $0$-connected graded braided bialgebra $\left( T=T(V,\mathfrak{c}%
),\nabla _{T},1_{T},\Delta _{T},\varepsilon _{T},\mathfrak{c}_{T}\right) $
where

\begin{itemize}
\item $\left( T=T(V,\mathfrak{c}),\nabla _{T},1_{T}\right) $ is the tensor
algebra $T(V)$ i.e. the free algebra generated by $V$.

\item $\mathfrak{c}_{T}$ is constructed iteratively from
$\mathfrak{c}$.

\item $\Delta _{T}:T\rightarrow T\otimes _{\mathfrak{c}_{T}}T$ is the unique
algebra homomorphism defined by setting $\Delta _{T}\left( v\right)
=1_{T}\otimes v+v\otimes 1_{T}$ for every $v\in V$.

\item $\varepsilon _{T}:T\rightarrow K$ is the unique algebra homomorphism
defined by setting $\varepsilon _{T}\left( v\right) =0$ for every $v\in V.$
\end{itemize}
\end{proposition}

\begin{remark}
\label{nn: T_c(V)}Note that $\Delta _{T}$ is the dual construction of the
quantum shuffle product, introduced by Rosso in \cite{Ro}. The $K$-linear
map $\mathfrak{c}_{T}^{n,m}$ is given in the figure below:
\begin{equation*}
\includegraphics{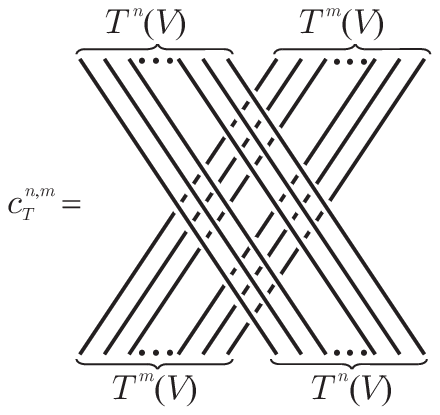}
\end{equation*}%
where each crossing represents a copy of $\mathfrak{c}$.
\end{remark}

\begin{remark}
\label{rem: 1.17}If $A=\bigoplus_{n\in \mathbb{N}}A^{n}$ is a $0$-connected
graded $\mathfrak{c}{}$-bialgebra then $\Delta ^{n,0}=\Delta ^{0,n}=\mathrm{%
Id}_{A^{n}}.$ Indeed, the proof given in the case $A=T(V)$ works for an
arbitrary connected graded $\mathfrak{c}{}$-bialgebra.
\end{remark}

\begin{theorem}
\label{teo: univ T(V,c)}Let $(V,\mathfrak{c}{})$ be a braided vector space.
Then $i_{V}:V\rightarrow T(V)$ is a morphism of braided vector spaces. If $%
(A,\nabla _{A},1_{A},\mathfrak{c}_{A})$ is a braided algebra and $%
f:V\rightarrow A$ is a morphism of braided vector spaces then there is a
unique morphism $\widetilde{f}:T(V,\mathfrak{c})\rightarrow A$ of braided
algebras that lifts $f$. If, in addition, $A$ is a braided bialgebra and $%
f(V)$ is contained in $P(A),$ the set of primitive elements of $A,$ then $%
\widetilde{f}$ is a morphism of braided bialgebras.
\end{theorem}

\begin{noname}
Recall (\cite[page 74]{Ka}) that the $X$-binomial coefficients $\binom{n}{k}%
_{X}$ are defined as follows. We set $(0)_{X}=(0)!_{X}=1$ and, for $n>0$,
define $(n)_{X}:=1+X+\cdots +X^{n-1}$ and $(n)!_{X}=(1)_{X}(2)_{X}\cdots
(n)_{X}$. Then:
\begin{equation}
\tbinom{n}{k}_{X}=\frac{(n)!_{X}}{(k)!_{X}\left( n-k\right) !_{X}}.
\label{ec:coeficienti binomiali}
\end{equation}%
It is well known that $\tbinom{n}{k}_{X}$ is a polynomial. Therefore we may
specialize $X$ at an arbitrary element $\lambda \in K.$ In this way we get
an element $\tbinom{n}{k}_{\lambda }\in K.$ Note that, if $\lambda $ is a
root a unity, then $\tbinom{n}{k}_{\lambda }$ may be zero. If $\mathrm{char}%
(K)=0$ and $\lambda =1,$ the formula shows that $\binom{n}{k}_{1}$ is the
classical binomial coefficient.
\end{noname}

\begin{definition}
\label{de:Hecke}We say that a braided vector space $(V,\mathfrak{c})$ is of
\emph{Hecke-type (or that }$c$\emph{\ is of Hecke-type) of mark} $\mathfrak{%
m(c)}=\lambda $ if
\begin{equation*}
(\mathfrak{c}+\mathrm{Id}_{V^{\otimes 2}})(\mathfrak{c}-\lambda \mathrm{Id}%
_{V^{\otimes 2}})=0.
\end{equation*}
\end{definition}

\begin{remark}
For every Hecke-type braiding $\mathfrak{c}$ of mark $\lambda ,$ the
operator $\check{\mathfrak{c}}:=-\lambda ^{-1}\mathfrak{c}$ is also of
Hecke-type. We have $\mathfrak{m}(\check{\mathfrak{c}})=\lambda ^{-1}$. Note
that, for $\mathfrak{d}=\check{\mathfrak{c}},$ one has $\check{\mathfrak{d}}=%
\mathfrak{c}.$
\end{remark}

\section{Braided enveloping algebras}

In this section we introduce the main example of $\mathfrak{c}$-bialgebras
that we will deal with, namely the enveloping algebra of a $\mathfrak{c}$%
-Lie algebra.

Given a vector spaces $V,W$ and a $K$-linear map $\alpha :V\otimes
V\rightarrow W$, we will denote $\alpha \otimes V$ and $V\otimes \alpha $ by
$\alpha _{1}$ and $\alpha _{2}$ respectively.

\begin{definition}
\label{def: braided bracket}Let $(V,\mathfrak{c})$ be a braided vector
space. We say that a $K$-linear map $\mathfrak{b}:V\otimes V\rightarrow V$
is a \emph{$\mathfrak{c}$-bracket}, or \emph{braided bracket}, if the
following compatibility conditions hold true:
\begin{equation}
\mathfrak{cb}_{1}=\mathfrak{b}_{2}\mathfrak{c}_{1}\mathfrak{c}_{2}\qquad
\text{and}\qquad \mathfrak{cb}_{2}=\mathfrak{b}_{1}\mathfrak{c}_{2}\mathfrak{%
c}_{1}  \label{ec: b compatibilitate}
\end{equation}%
Let $\mathfrak{b}$ be a $\mathfrak{c}$-bracket on $(V,\mathfrak{c})$ and let
$\mathfrak{b}^{\prime }$ be a $\mathfrak{c}^{\prime }$-bracket on $%
(V^{\prime },\mathfrak{c}^{\prime })$. We will say that a morphism of
braided vector spaces from $V$ to $V^{\prime }$ is \emph{a morphism of
braided brackets} if $f\mathfrak{b}=\mathfrak{b}^{\prime }(f\otimes f)$.

Let $\mathfrak{b}$ be a bracket on a braided vector space $\left( V,%
\mathfrak{c}\right) $ of Hecke-type with mark $\lambda $. Let $I_{\mathfrak{%
c,b{}}}$ is the two-sided ideal generated by the set
\begin{equation*}
X_{\mathfrak{c},\mathfrak{b}}=\{\mathfrak{c}(z)-\lambda z-\mathfrak{b}%
(z)\mid z\in V\otimes V=T^{2}(V)\}.
\end{equation*}%
The \emph{enveloping algebra of} $(V,\mathfrak{c},\mathfrak{b})$ is by
definition the algebra%
\begin{equation*}
U=U(V,\mathfrak{c},\mathfrak{b}):=\frac{T(V,{\mathfrak{c}})}{I_{\mathfrak{c},%
\mathfrak{b}}}.
\end{equation*}%
We will denote by $\pi _{\mathfrak{c,b}}:T(V,{\mathfrak{c}})\rightarrow U(V,%
\mathfrak{c},\mathfrak{b})$ the canonical projection.

When $\mathfrak{b}=0,$ the enveloping algebra of $(V,\mathfrak{c},0)$ is
called the $\mathfrak{c}$-\emph{symmetric algebra}, or \emph{braided
symmetric algebra}, if there is no danger of confusion. It will be denoted
by $S(V,{\mathfrak{c}})$.
\end{definition}

\begin{proposition}
\label{teo: U(V,c,b) bialg}Let $\mathfrak{b}$ be a bracket on a braided
vector space $\left( V,\mathfrak{c}\right) $ of Hecke-type with mark $%
\lambda .$ Then $I_{\mathfrak{c},\mathfrak{b{}}}$ is a coideal in $T(V,{%
\mathfrak{c}})$. Moreover, on the quotient algebra (and coalgebra) there is
a natural braiding $\mathfrak{c}_{U}$ such that $(U(V,\mathfrak{c},\mathfrak{%
b}),\mathfrak{c}_{U})$ is a braided bialgebra.
\end{proposition}

\begin{proof}
We denote $T(V,\mathfrak{c})$ and $T^{n}(V,\mathfrak{c})$ by $T$ and $T^{n}$%
, respectively. We first prove that ${\mathfrak{c}}_{T}$ maps $I_{\mathfrak{c%
},\mathfrak{b}}\otimes T$ into $T\otimes I_{\mathfrak{c},\mathfrak{b}}$ and $%
T\otimes I_{\mathfrak{c},\mathfrak{b}}$ into $I_{\mathfrak{c},\mathfrak{b}%
}\otimes T.$ Let $x\in T^{n},$ $y\in T^{m},$ $t\in T^{p}$ and $z\in T^{2}.$
Since $\mathfrak{c}$ verifies the braid equation and (\ref{ec: b
compatibilitate}) we get:
\begin{eqnarray*}
{\mathfrak{c}}^{n+m+2,p}(T^{n}\otimes \mathfrak{c}\otimes T^{m}\otimes
T^{p}) &=&(T^{p}\otimes T^{n}\otimes \mathfrak{c}\otimes T^{m}){\mathfrak{c}}%
^{n+m+2,p}, \\
{\mathfrak{c}}^{n+m+1,p}(T^{n}\otimes \mathfrak{b}\otimes T^{m}\otimes
T^{p}) &=&(T^{p}\otimes T^{n}\otimes \mathfrak{b}\otimes T^{m}){\mathfrak{c}}%
^{n+m+2,p}.
\end{eqnarray*}%
Then%
\begin{equation*}
\mathfrak{c}_{T}\left[ \left( T^{n}\otimes (\mathfrak{c}-\lambda \mathrm{Id}%
_{T^{2}}-\mathfrak{b})\otimes T^{m}\right) \otimes T^{p}\right] =\left[
T^{p}\otimes \left( T^{n}\otimes (\mathfrak{c}-\lambda \mathrm{Id}_{T^{2}}-%
\mathfrak{b})\otimes T^{m}\right) \right] \mathfrak{c}_{T},
\end{equation*}%
relation that shows us that ${\mathfrak{c}}_{T}$ maps $I_{\mathfrak{c},%
\mathfrak{b{}}}\otimes T$ into $T\otimes I_{\mathfrak{c},\mathfrak{b{}}}$.
The other property can be proved similarly.

We claim that $X_{\mathfrak{c},\mathfrak{b}}$ is a coideal in $T$. In fact
we will prove that $X_{\mathfrak{c},\mathfrak{b}}$ contains only primitive
elements in $T$. Let $z\in T^{2}$. By Proposition \ref{pr: Delta(c)} we have
\begin{eqnarray*}
\Delta _{T}(z) &=&z\otimes 1+z+\mathfrak{c}(z)+1\otimes z, \\
\Delta _{T}(\mathfrak{c}(z)) &=&\mathfrak{c}(z)\otimes 1+\mathfrak{c}(z)+%
\mathfrak{c}^{2}(z)+1\otimes \mathfrak{c}(z).
\end{eqnarray*}%
Thus $\Delta _{T}(\mathfrak{c}(z)-\lambda z)=(\mathfrak{c}(z)-\lambda
z)\otimes 1+1\otimes (\mathfrak{c}(z)-\lambda z),$ as $\mathfrak{c}$ is a
Hecke braiding of mark $\lambda .$ It follows that $\mathfrak{c}(z)-\lambda
z\in P(T)$. Since $\mathfrak{b}(z)\in V$ we deduce that $\mathfrak{c}%
(z)-\lambda z-\mathfrak{b}(z)\in P(T)$, so $X_{\mathfrak{c},\mathfrak{b}%
}\subseteq P(T)$.

Now, by (\ref{ec: Produs1}), it easily follows that $I_{\mathfrak{c},%
\mathfrak{b}}$, the ideal generated by $X_{\mathfrak{c},\mathfrak{b}}$, is a
coideal. It remains to show that $\mathfrak{c}$ factors through a braiding $%
\mathfrak{c}_{U}$ of $U:=U(V,\mathfrak{c},\mathfrak{b})$, that makes $U$ a
braided bialgebra. Let $\pi _{\mathfrak{c,b}}:T\rightarrow U$ be the
canonical projection. By the foregoing, ${\mathfrak{c}}_{T}$ maps the kernel
of $\pi _{\mathfrak{c},\mathfrak{b}}\otimes \pi _{\mathfrak{c},b}$ into
itself, so there is a $K$-linear morphism $\mathfrak{c}_{U}:U\otimes
U\rightarrow U\otimes U$ such that
\begin{equation*}
\mathfrak{c}_{U}(\pi _{\mathfrak{c},\mathfrak{b}}\otimes \pi _{\mathfrak{c},%
\mathfrak{b{}}})=(\pi _{\mathfrak{c},\mathfrak{b}}\otimes \pi _{\mathfrak{c},%
\mathfrak{b}}){\mathfrak{c}}_{T}.
\end{equation*}%
Since $T$ is a ${\mathfrak{c}}_{T}$-bialgebra, this relation entails that $U$
is a $\mathfrak{c}_{U}$-bialgebra and that the canonical projection $\pi _{%
\mathfrak{c},\mathfrak{b}}$ becomes a morphism of braided bialgebras.
\end{proof}

\begin{remark}
Let $(V,\mathfrak{c})$ be an arbitrary braided vector space. Let $K^{n}:=%
\mathrm{Ker}(\mathfrak{S}_{n})$ where $\mathfrak{S}_{n}$ denotes the quantum
symmetrizer \cite[2.3]{AS2}. It is well-known that $\bigoplus_{n\mathbb{N}%
}K^{n}$ is an ideal and a coideal in $T(V,\mathfrak{c}),$ see \cite[section 3%
]{AS2}$.$ Since
\begin{equation*}
\mathfrak{c}_{T}^{n,m}(T^{n}\otimes K^{m}+K^{n}\otimes T^{m})\subseteq
K^{m}\otimes T^{n}+T^{m}\otimes K^{n}
\end{equation*}%
it follows that $B(V,\mathfrak{c})=T(V,\mathfrak{c})/\left( \bigoplus_{n\in
\mathbb{N}}K^{n}\right) $ is a quotient graded braided bialgebra of $T(V,%
\mathfrak{c})$, that is called the Nichols algebra of $(V,\mathfrak{c})$.
Let us now assume that $\mathfrak{c}{}$is a braiding of Hecke-type of mark $%
\lambda .$ By the definition of Hecke operators we have
\begin{equation*}
\mathrm{Im}(\mathfrak{c}-\lambda \mathrm{Id}_{T^{2}})\subseteq \mathrm{Ker}(%
\mathrm{Id}_{T^{2}}+\mathfrak{c})=K^{2}\subseteq \bigoplus_{n\in \mathbb{N}%
}K^{n}.
\end{equation*}%
Therefore, there is a morphism of braided bialgebras $\varphi :S(V,\mathfrak{%
c})\rightarrow B(V,\mathfrak{c})$ such that $\varphi |_{V}=\mathrm{Id}_{V}.$
Obviously $\varphi $ is surjective, since $B(V,\mathfrak{c})$ is generated
by $V.$ Later (see Theorem \ref{teo: S(V,c) type one}) we will see that the
space of primitive elements in $S(V,\mathfrak{c})$ and the homogeneous
component $S^{1}(V,\mathfrak{c})=V$ are identical. By \cite[Theorem 5.3.1]%
{Mo}, it follows that $\varphi $ is injective too. Thus, $S(V,\mathfrak{c})$
and $B(V,\mathfrak{c})$ are isomorphic braided bialgebras.

We are going to investigate some basic properties of these objects.
\end{remark}

\begin{proposition}
\label{pr: universalitate U}Let $(A,\nabla ,1,\Delta ,\varepsilon ,\mathfrak{%
c}_{A})$ be a connected braided bialgebra. Let $P$ be the space of primitive
elements of $A.$ Assume that there is $\lambda \in K^{\ast }$ such that $%
\mathfrak{c}_{P}:=\mathfrak{c}_{A}|_{P\otimes P}$ is a braiding of
Hecke-type on $P$ of mark $\lambda .$ Then:

a) $\nabla (\mathfrak{c}_{P}-\lambda \mathrm{Id}_{P^{\otimes 2}})(P\otimes
P)\subseteq P,$ so we can define $\mathfrak{b}_{P}:P\otimes P\rightarrow P$
by $\mathfrak{b}_{P}=\nabla (\mathfrak{c}_{P}-\lambda \mathrm{Id}%
_{P^{\otimes 2}})|_{P\otimes P}.$

b) The map $\mathfrak{b}_{P}$ is a braided bracket on the braided vector
space $(P,\mathfrak{c}_{P})$.

c) Let $f:(V,\mathfrak{c},\mathfrak{b})\rightarrow (P,\mathfrak{c}_{P},%
\mathfrak{b}_{P})$ be a morphism of braided brackets and assume that $%
\mathfrak{c}$ is a braiding of Hecke-type with mark $\lambda .$ Then there
is a unique morphism of braided bialgebras $\widetilde{f}:U(V,\mathfrak{c},%
\mathfrak{b})\rightarrow A$ that lifts $f$.
\end{proposition}

\begin{proof}
First, observe that, by Remark \ref{rem: c_P}, $\mathfrak{c}_{A}(P\otimes
P)\subseteq P\otimes P$ so that it makes sense to consider $\mathfrak{c}%
_{P}:P\otimes P\rightarrow P\otimes P.$

a) By assumption, $\mathfrak{c}_{P}$ is a Hecke operator on $P$ of mark $%
\lambda .$ For $z\in P\otimes P,$ by (\ref{ec: Produs1}) we get:%
\begin{eqnarray*}
\Delta \nabla (z) &=&\nabla \left( z\right) \otimes 1+z+\mathfrak{c}%
_{P}(z)+1\otimes \nabla \left( z\right) , \\
\Delta \nabla \mathfrak{c}_{P}(z) &=&\nabla \mathfrak{c}_{P}\left( z\right)
\otimes 1+\mathfrak{c}_{P}\left( z\right) +\mathfrak{c}_{P}^{2}(z)+1\otimes
\nabla \mathfrak{c}_{P}\left( z\right) .
\end{eqnarray*}%
Therefore $\nabla (\mathfrak{c}_{P}-\lambda \mathrm{Id}_{P^{\otimes
2}})\left( z\right) \in P.$ This shows that $\nabla (\mathfrak{c}%
_{P}-\lambda \mathrm{Id}_{P^{\otimes 2}})(P\otimes P)\subseteq P.$

b) We have to prove the compatibility relation between $\mathfrak{c}$ and $%
\mathfrak{b{}},$ that is we have (\ref{ec: b compatibilitate}). But these
relations follows immediately by the braid relation and the fact that $A$ is
a $\mathfrak{c}_{A}$-algebra.

c) Apply the universal property of $T(V,\mathfrak{c})$ (see Theorem \ref%
{teo: univ T(V,c)}) to get a morphism $f^{\prime }:T(V,\mathfrak{c}%
)\rightarrow A$ of braided bialgebras that lifts $f$. Since $f$ is a
morphism of braided brackets, by the definition of $\mathfrak{b}_{P}$, it
results that $f^{\prime }$ maps $\mathfrak{c}(z)-\lambda z-\mathfrak{b}(z)$
to $0$. Therefore $f^{\prime }$ factors through a morphism $\widetilde{f}%
:U(V,\mathfrak{c},\mathfrak{b})\rightarrow A$, which lifts $f$ and is
compatible with the braidings (note that $U(V,\mathfrak{c},\mathfrak{b})$ is
a braided bialgebra in view of Theorem \ref{teo: U(V,c,b) bialg}).
\end{proof}

\begin{proposition}
\label{pr: universalitate S}Let $(A,\nabla ,1,\Delta ,\varepsilon ,\mathfrak{%
c}{}_{A})$ be a connected braided bialgebra. Let $P$ denote the primitive
part of $A$. Assume that there is $\lambda \in K^{\ast }$ such that $\nabla
\mathfrak{c}_{A}=\lambda \nabla $ on $P\otimes P$. If $\mathfrak{c}_{P}=%
\mathfrak{c}_{A}|_{P\otimes P}$ then $\mathfrak{c}_{P}$ is of Hecke-type on $%
P$ of mark $\lambda $. Moreover, if $(V,\mathfrak{c})$ is a braided vector
space such that $\mathfrak{c}$ is a Hecke operator of mark $\lambda $ and $%
f:(V,\mathfrak{c})\rightarrow (P,\mathfrak{c}_{P})$ is a morphism of braided
vector spaces then there is a unique morphism of braided bialgebras $%
\widetilde{f}:S(V,\mathfrak{c})\rightarrow A$ that lifts $f$. If $A$ is
graded then $\widetilde{f}$ respects the gradings on $S(V,\mathfrak{c})$ and
$A.$
\end{proposition}

\begin{proof}
For $z\in P\otimes P,$ by (\ref{ec: Produs1}) we get:%
\begin{eqnarray*}
\Delta \nabla (z) &=&\nabla \left( z\right) \otimes 1+z+\mathfrak{c}%
_{A}(z)+1\otimes \nabla \left( z\right) , \\
\Delta \nabla \mathfrak{c}_{A}(z) &=&\nabla \mathfrak{c}_{A}\left( z\right)
\otimes 1+\mathfrak{c}_{A}\left( z\right) +\mathfrak{c}_{A}^{2}(z)+1\otimes
\nabla \mathfrak{c}_{A}\left( z\right) .
\end{eqnarray*}%
By assumption, $\nabla \mathfrak{c}_{A}(z)=\lambda \nabla (z)$ whence
\begin{equation*}
0=\mathfrak{c}_{A}\left( z\right) +\mathfrak{c}_{A}^{2}(z)-\lambda z-\lambda
\mathfrak{c}_{A}(z)=\left( \mathfrak{c}_{A}+\mathrm{Id}_{A^{\otimes
2}}\right) \left( \mathfrak{c}_{A}-\lambda \mathrm{Id}_{A^{\otimes
2}}\right) (z)
\end{equation*}%
Thus $\mathfrak{c}_{P}$ is a Hecke operator of mark $\lambda $. By taking $%
\mathfrak{b}=0$ in Proposition \ref{pr: universalitate U}, it results that
there is $\widetilde{f}$ that lifts $f$.
\end{proof}

\begin{remark}
\label{re: universalitate S}The above proposition still works under the
slighter assumption $\nabla \mathfrak{c}_{A}=\lambda \nabla $ on $\mathrm{Im}%
\,f\otimes \mathrm{Im}\,f$.
\end{remark}

\begin{noname}
Let $T^{n}:=T^{n}(V)$ and let $T^{\leq n}:=\bigoplus_{0\leq m\leq n}T^{m}.$
By construction $\pi _{\mathfrak{c},\mathfrak{b}}$ is a morphism of algebras
and coalgebras from $T(V,\mathfrak{c})$ to $U(V,\mathfrak{c},\mathfrak{b})$.
Thus $U_{n}^{\prime }:=\pi _{\mathfrak{c},\mathfrak{b}}(T^{\leq n})$ defines
a braided bialgebra filtration on $U(V,\mathfrak{c},\mathfrak{b})$, i.e. $%
(U_{n}^{\prime })_{n\in \mathbb{N}}$ is an algebra and coalgebra filtration
on $U(V,\mathfrak{c},\mathfrak{b})$ which is compatible with $\mathfrak{c}%
_{U}$. It will be called the \emph{standard filtration }on $U(V,\mathfrak{c},%
\mathfrak{b}).$ In general, this filtration and the coradical filtration \ $%
(U_{n})_{n\in \mathbb{N}}$ are not identical, but we always have $%
U_{n}^{\prime }\subseteq U_{n},$ for any $n\in \mathbb{N}.$

If $\mathfrak{b}=0$ then $S(V,\mathfrak{c}):=U(V,\mathfrak{c},0)$ is a
graded $\mathfrak{c}_{S}$-bialgebra, $S(V,\mathfrak{c})=\bigoplus_{n\in
\mathbb{N}}S^{n}(V,\mathfrak{c})$. The standard filtration on $S(V,\mathfrak{%
c})$ is the filtration associated to this grading.
\end{noname}

\begin{proposition}
\label{pr: Sconn}Let $\mathfrak{b}$ be a $\mathfrak{c}$-bracket on a braided
vector space $(V,\mathfrak{c})$ of Hecke-type. Then $U(V,\mathfrak{c},%
\mathfrak{b})$ is a connected coalgebra. Moreover, for every braided vector
space $(V,\mathfrak{c}),$ $S(V,\mathfrak{c})$ is a $0$-connected graded
braided coalgebra.
\end{proposition}

\begin{proof}
We know that the tensor algebra of an arbitrary braided vector space is a $0$%
-connected coalgebra. By definition, $U(V,\mathfrak{c},\mathfrak{b})$ is a
quotient coalgebra of $T(V,\mathfrak{c}),$ where $\lambda =\mathfrak{m}(%
\mathfrak{c}).$ Then, in view of \cite[Corollary 5.3.5]{Mo}, $U(V,\mathfrak{c%
},\mathfrak{b})$ is connected. In particular, braided symmetric algebras are
connected coalgebras. They are also $0$-connected since they are graded
quotients of $T(V,\mathfrak{c})$.
\end{proof}

\begin{remark}
Let $\mathfrak{b}$ be a $\mathfrak{c}$-bracket on a braided vector space $(V,%
\mathfrak{c})$ of Hecke-type. The composition of the inclusion $V\rightarrow
T(V,\mathfrak{c})$ with the canonical projection $\pi _{\mathfrak{c}{},%
\mathfrak{b{}}}$ gives a map
\begin{equation*}
\iota _{\mathfrak{c},\mathfrak{b}}:V\rightarrow U(V,\mathfrak{c},\mathfrak{b}%
).
\end{equation*}%
Its image is included in the space of primitive elements of $U(V,\mathfrak{c}%
,\mathfrak{b}).$ In general $\iota _{\mathfrak{c}{},\mathfrak{b}}$ is
neither injective nor onto. Our purpose now is to investigate when $\iota _{%
\mathfrak{c},\mathfrak{b}}$ is injective (see Theorem (\ref{te: Lie algebras}%
)).
\end{remark}

\begin{definition}
Let $\left( B,\nabla _{B},1_{B},\Delta _{B},\varepsilon _{B},c_{B}\right) $
be a graded braided bialgebra. For every $a,b,n\in
%TCIMACRO{\U{2115} }%
%BeginExpansion
\mathbb{N}
%EndExpansion
$, set
\begin{equation*}
\mathrm{\Gamma }_{a,b}^{B}:=\nabla _{B}^{a,b}\Delta _{B}^{a,b}\text{.}
\end{equation*}
\end{definition}

\begin{lemma}
Let $\left( B,\nabla _{B},1_{B},\Delta _{B},\varepsilon _{B},\mathfrak{c}%
_{B}\right) $ be a $0$-connected graded braided bialgebra. Then%
\begin{equation*}
\Delta _{B}^{n,0}\left( z\right) =z\otimes 1_{B}\qquad \text{and}\qquad
\Delta _{B}^{0,n}=1_{B}\otimes z,\text{ for every }z\in B^{n}.
\end{equation*}%
Moreover%
\begin{equation}
\Delta _{B}^{n,1}\nabla _{B}^{n,1}=\mathrm{Id}_{B^{n}\otimes B^{1}}+\left(
\nabla _{B}^{n-1,1}\otimes B^{1}\right) \left( B^{n-1}\otimes \mathfrak{c}%
_{B}^{1,1}\right) \left( \Delta _{B}^{n-1,1}\otimes B^{1}\right) .
\label{form: Delta Mult}
\end{equation}
\end{lemma}

\begin{proof}
The first assertion follows by Remark \ref{rem: 1.17}. Since $B$ is a graded
bialgebra we have%
\begin{eqnarray*}
&&\Delta _{B}^{n,1}\nabla _{B}^{n,1}=\left[
\begin{array}{c}
\left( \nabla _{B}^{n,0}\otimes \nabla _{B}^{0,1}\right) \left( B^{n}\otimes
\mathfrak{c}_{B}^{0,0}\otimes B^{1}\right) \left( \Delta _{B}^{n,0}\otimes
\Delta _{B}^{0,1}\right) + \\
+\left( \nabla _{B}^{n-1,1}\otimes \nabla _{B}^{1,0}\right) \left(
B^{n-1}\otimes \mathfrak{c}_{B}^{1,1}\otimes B^{0}\right) \left( \Delta
_{B}^{n-1,1}\otimes \Delta _{B}^{1,0}\right)%
\end{array}%
\right] \\
&=&\mathrm{Id}_{B^{n}\otimes B^{1}}+\left( \nabla _{B}^{n-1,1}\otimes
B^{1}\right) \left( B^{n-1}\otimes \mathfrak{c}_{B}^{1,1}\right) \left(
\Delta _{B}^{n-1,1}\otimes B^{1}\right) ,
\end{eqnarray*}%
so that (\ref{form: Delta Mult}) holds.
\end{proof}

\begin{definition}
\label{def: strongly grAlg}Let $(A,\nabla ,1_{A})$ be a graded algebra. We
say that $A$ is a \emph{strongly }$%
%TCIMACRO{\U{2115} }%
%BeginExpansion
\mathbb{N}
%EndExpansion
$\emph{-graded algebra} whenever $\nabla _{i,j}:A_{i}\otimes
A_{j}\rightarrow A_{i+j}$ is an epimorphism for every $i,j\in \mathbb{N}$
(equivalently $\nabla _{n,1}:A_{n}\otimes A_{1}\rightarrow A_{n+1}$ is an
epimorphism for every $n\in \mathbb{N}$).

Dually, let $(C,\Delta ,\varepsilon )$ be a graded coalgebra. We say that $C$
is a \emph{strongly }$%
%TCIMACRO{\U{2115} }%
%BeginExpansion
\mathbb{N}
%EndExpansion
$\emph{-graded coalgebra} whenever $\Delta _{i,j}:C_{i+j}\rightarrow
C_{i}\otimes C_{j}$ is a monomorphism for every $i,j\in \mathbb{N}$
(equivalently $\Delta _{n,1}:C_{n+1}\rightarrow C_{n}\otimes C_{1}$ is a
monomorphism for every $n\in \mathbb{N}$).

For more details on these (co)algebras see e.g. \cite{AM- type1}.
\end{definition}

\begin{definition}
\label{def: regular}An element $\lambda \in K^{\ast }$ is called $n$-\textbf{%
regular} whenever $(k)_{\lambda }\neq 0$, for any $1\leq k\leq n$. If $%
\lambda $ is $n$-regular for any $n>0,$ we will simply say that $\lambda $
is \textbf{regular}.
\end{definition}

\begin{remark}
If $\lambda $ is $n$-regular (respectively regular) then $\lambda ^{-1}$ is
also $n$-regular (respectively regular).
\end{remark}

\begin{theorem}
\label{teo: aureo}Let $\left( B,\nabla _{B},1_{B},\Delta _{B},\varepsilon
_{B},\mathfrak{c}_{B}\right) $ be a $0$-connected graded braided bialgebra
and let $\lambda \in K^{\ast }$ be regular. The following are equivalent.

\begin{enumerate}
\item $B$ is a bialgebra of type one and $\mathfrak{c}_{B}^{1,1}$ is a
braiding of Hecke-type of mark $\lambda $.

\item $B$ is strongly $%
%TCIMACRO{\U{2115} }%
%BeginExpansion
\mathbb{N}
%EndExpansion
$-graded as a coalgebra, $\nabla _{B}^{1,1}$ is surjective and $\mathfrak{c}%
_{B}^{1,1}$ is a braiding of Hecke-type of mark $\lambda $.

\item $B$ is strongly $%
%TCIMACRO{\U{2115} }%
%BeginExpansion
\mathbb{N}
%EndExpansion
$-graded as a coalgebra and $\left( \mathfrak{c}_{B}^{1,1}-\lambda \mathrm{Id%
}_{B^{2}}\right) \Delta _{B}^{1,1}=0$.

\item $B$ is strongly $%
%TCIMACRO{\U{2115} }%
%BeginExpansion
\mathbb{N}
%EndExpansion
$-graded as an algebra, $\Delta _{B}^{1,1}$ is injective and $\mathfrak{c}%
_{B}^{1,1}$ is a braiding of Hecke-type of mark $\lambda $.

\item $B$ is strongly $%
%TCIMACRO{\U{2115} }%
%BeginExpansion
\mathbb{N}
%EndExpansion
$-graded as an algebra and $\nabla _{B}^{1,1}\left( \mathfrak{c}%
_{B}^{1,1}-\lambda \mathrm{Id}_{B^{2}}\right) =0$.
\end{enumerate}
\end{theorem}

\begin{proof}
$\left( 1\right) \Rightarrow \left( 2\right) $ By definition, $B$ is
strongly $%
%TCIMACRO{\U{2115} }%
%BeginExpansion
\mathbb{N}
%EndExpansion
$-graded both as a coalgebra and as an algebra.

$\left( 2\right) \Rightarrow \left( 3\right) $ We have
\begin{equation*}
0=\left( \mathfrak{c}_{B}^{1,1}-\lambda \mathrm{Id}_{B^{2}}\right) \left(
\mathfrak{c}_{B}^{1,1}+\mathrm{Id}_{B^{2}}\right) \overset{\text{(\ref{form:
Delta Mult})}}{=}\left( \mathfrak{c}_{B}^{1,1}-\lambda \mathrm{Id}%
_{B^{2}}\right) \Delta _{B}^{1,1}\nabla _{B}^{1,1}.
\end{equation*}%
Since $\nabla _{B}^{1,1}$ is an epimorphism, we get $\left( \mathfrak{c}%
_{B}^{1,1}-\lambda \mathrm{Id}_{B^{2}}\right) \Delta _{B}^{1,1}=0$.

$\left( 3\right) \Rightarrow \left( 1\right) $ We have%
\begin{eqnarray*}
&&\Delta _{B}^{n,1}\mathrm{\Gamma }_{n,1}^{B}=\Delta _{B}^{n,1}\nabla
_{B}^{n,1}\Delta _{B}^{n,1}\overset{\text{(\ref{form: Delta Mult})}}{=}%
\Delta _{B}^{n,1}+\left( \nabla _{B}^{n-1,1}\otimes B^{1}\right) \left(
B^{n-1}\otimes \mathfrak{c}_{B}^{1,1}\right) \left( \Delta
_{B}^{n-1,1}\otimes B^{1}\right) \Delta _{B}^{n,1} \\
&=&\Delta _{B}^{n,1}+\left( \nabla _{B}^{n-1,1}\otimes B^{1}\right) \left(
B^{n-1}\otimes \mathfrak{c}_{B}^{1,1}\Delta _{B}^{1,1}\right) \Delta
_{B}^{n-1,2} \\
&=&\Delta _{B}^{n,1}+\left( \nabla _{B}^{n-1,1}\otimes B^{1}\right) \left(
B^{n-1}\otimes \lambda \Delta _{B}^{1,1}\right) \Delta _{B}^{n-1,2} \\
&=&\Delta _{B}^{n,1}+\lambda \left( \nabla _{B}^{n-1,1}\Delta
_{B}^{n-1,1}\otimes B^{1}\right) \Delta _{B}^{n,1}=\Delta _{B}^{n,1}+\lambda
\left( \mathrm{\Gamma }_{n-1,1}^{B}\otimes B^{1}\right) \Delta _{B}^{n,1}
\end{eqnarray*}%
so that%
\begin{equation}
\Delta _{B}^{n,1}\mathrm{\Gamma }_{n,1}^{B}=\Delta _{B}^{n,1}+\lambda \left(
\mathrm{\Gamma }_{n-1,1}^{B}\otimes B^{1}\right) \Delta _{B}^{n,1}
\label{form: Gamma a}
\end{equation}%
Let us prove by induction that
\begin{equation}
\mathrm{\Gamma }_{n,1}^{B}=\left( n+1\right) _{\lambda }\mathrm{Id}%
_{B^{n+1}},\text{ for every }n\geq 1\text{.}  \label{form: gamma}
\end{equation}%
$n=1)$ We have
\begin{equation*}
\Delta _{B}^{1,1}\mathrm{\Gamma }_{1,1}^{B}\overset{\text{(\ref{form: Gamma
a})}}{=}\Delta _{B}^{1,1}+\lambda \Delta _{B}^{1,1}=\left( 2\right)
_{\lambda }\Delta _{B}^{1,1}.
\end{equation*}%
Since, by hypothesis, $\Delta _{B}^{1,1}$ is injective, we obtain $\mathrm{%
\Gamma }_{1,1}^{B}=\left( 2\right) _{q}\mathrm{Id}_{B^{2}}.$

$n-1\Rightarrow n)$ We have%
\begin{equation*}
\Delta _{B}^{n,1}\mathrm{\Gamma }_{n,1}^{B}\overset{\text{(\ref{form: Gamma
a})}}{=}\Delta _{B}^{n,1}+\lambda \left( \mathrm{\Gamma }_{n-1,1}^{B}\otimes
B^{1}\right) \Delta _{B}^{n,1}=\Delta _{B}^{n,1}+\lambda \left( n\right)
_{\lambda }\Delta _{B}^{n,1}=\left( n+1\right) _{\lambda }\Delta _{B}^{n,1}
\end{equation*}%
Since, by hypothesis, $\Delta _{B}^{n,1}$ is injective, we obtain $\mathrm{%
\Gamma }_{n,1}^{B}=\left( n+1\right) _{\lambda }\mathrm{Id}_{B^{n+1}}.$

We have so proved (\ref{form: gamma}).

Since $\lambda $ is regular, we have $\left( n+1\right) _{\lambda }\neq 0$
so that
\begin{equation*}
\nabla _{B}^{n,1}\Delta _{B}^{n,1}=\mathrm{\Gamma }_{n,1}^{B}\overset{\text{(%
\ref{form: gamma})}}{=}\left( n+1\right) _{\lambda }\mathrm{Id}_{B^{n+1}}
\end{equation*}%
is bijective for every $n\geq 1.$ Therefore $\nabla _{B}^{n,1}$ is
surjective for every $n\geq 1.$ Equivalently $B$ is strongly $%
%TCIMACRO{\U{2115} }%
%BeginExpansion
\mathbb{N}
%EndExpansion
$-graded as an algebra and hence of type one. Moreover
\begin{equation*}
\left( \mathfrak{c}_{B}^{1,1}-\lambda \mathrm{Id}_{B^{2}}\right) \left(
\mathfrak{c}_{B}^{1,1}+\mathrm{Id}_{B^{2}}\right) \overset{\text{(\ref{form:
Delta Mult})}}{=}\left( \mathfrak{c}_{B}^{1,1}-\lambda \mathrm{Id}%
_{B^{2}}\right) \Delta _{B}^{1,1}\nabla _{B}^{1,1}=0.
\end{equation*}

$\left( 1\right) \Leftrightarrow \left( 4\right) \Leftrightarrow \left(
5\right) $ It follows by dual arguments.
\end{proof}

\begin{definition}
A graded coalgebra $C=\bigoplus_{n\in \mathbb{N}}C^{n}$ is called \textbf{%
strict} if it is $0$-connected and $P(C)=C^{1}.$
\end{definition}

\begin{theorem}
\label{teo: S(V,c) type one}(cf. \cite[Proposition 3.4]{AS2}) Let $\left( V,%
\mathfrak{c}\right) $ be a braided vector space of Hecke-type with regular
mark $\lambda .$ Then $S\left( V,\mathfrak{c}\right) $ is a bialgebra of
type one. In particular $S(V,\mathfrak{c})$ is a strict coalgebra.
\end{theorem}

\begin{proof}
By definition, we have
\begin{equation*}
S:=S\left( V,\mathfrak{c}\right) =U(V,\mathfrak{c},0)=\frac{T(V,{\mathfrak{c}%
})}{\left( \mathfrak{c}(z)-\lambda z\mid z\in V\otimes V\right) }.
\end{equation*}%
Thus, since $S$ is a graded quotient of the graded braided bialgebra $%
T\left( V,\mathfrak{c}\right) ,$ we get that $S$ is strongly $%
%TCIMACRO{\U{2115} }%
%BeginExpansion
\mathbb{N}
%EndExpansion
$-graded as an algebra. Moreover $\nabla _{S}^{1,1}\left( \mathfrak{c}%
_{S}^{1,1}-\lambda \mathrm{Id}_{S^{2}}\right) =0.$ By Theorem \ref{teo:
aureo}, we conclude.
\end{proof}

\section{Categorical subspaces}

\begin{definition}
\cite[2.2]{Kharchenko- connected} A subspace $L$ of a braided vector space $%
\left( V,\mathfrak{c}\right) $ is said to be \textbf{categorical} if
\begin{equation}
\mathfrak{c}\left( L\otimes V\right) \subseteq V\otimes L\qquad \text{and}%
\qquad \mathfrak{c}\left( V\otimes L\right) \subseteq L\otimes V.
\label{form: categorical}
\end{equation}
\end{definition}

\begin{theorem}
\label{teo: categorical}Let $\left( V,\mathfrak{c}\right) $ be a braided
vector space of Hecke-type of mark $\lambda $. Assume $\lambda \neq 0,1.$
Let $L$ be a categorical subspace of $V$. Then $L=0$ or $L=V$.
\end{theorem}

\begin{proof}
From $\left( \mathfrak{c}+\mathrm{Id}_{V\otimes V}\right) \left( \mathfrak{c}%
-\lambda \mathrm{Id}_{V\otimes V}\right) =0$ and $\lambda \neq 1,$ we get $%
\mathfrak{c}=\left( \lambda -1\right) ^{-1}\left( \mathfrak{c}^{2}-\lambda
\mathrm{Id}_{V\otimes V}\right) $ so that
\begin{equation*}
\mathfrak{c}\left( L\otimes V\right) =\frac{1}{\lambda -1}\left( \mathfrak{c}%
^{2}-\lambda \mathrm{Id}_{V\otimes V}\right) \left( L\otimes V\right)
\overset{\text{(\ref{form: categorical})}}{\subseteq }L\otimes V.
\end{equation*}
Then%
\begin{equation}
\mathfrak{c}\left( L\otimes V\right) \subseteq \left( V\otimes L\right) \cap
\left( L\otimes V\right) =L\otimes L.  \label{form: categorical2}
\end{equation}%
From $\left( \mathfrak{c}+\mathrm{Id}_{V\otimes V}\right) \left( \mathfrak{c}%
-\lambda \mathrm{Id}_{V\otimes V}\right) =0$ and $\lambda \neq 0,$ we get $%
\mathrm{Id}_{V\otimes V}=\lambda ^{-1}\left( \mathfrak{c}^{2}+\left(
1-\lambda \right) \mathfrak{c}\right) $ so that
\begin{equation*}
L\otimes V=\frac{1}{\lambda }\left( \mathfrak{c}^{2}+\left( 1-\lambda
\right) \mathfrak{c}\right) \left( L\otimes V\right) \overset{\text{(\ref%
{form: categorical2}),(\ref{form: categorical})}}{\subseteq }L\otimes L.
\end{equation*}
Since $L\subseteq V,$ we deduce that $L=0$ or $L=V.$
\end{proof}

\begin{proposition}
\label{pro: triviality}Let $\mathfrak{b}$ be a $\mathfrak{c}$-bracket on a
braided vector space $(V,\mathfrak{c})$ of Hecke-type of mark $\lambda $.
Assume $\lambda \neq 0,1.$ Then $\mathfrak{b}$ is zero or surjective.
\end{proposition}

\begin{proof}
Let $L=\text{Im}\left( \mathfrak{b}\right) \subseteq V.$ We have that%
\begin{eqnarray*}
\mathfrak{cb}_{1} &=&\mathfrak{b}_{2}\mathfrak{c}_{1}\mathfrak{c}%
_{2}\Rightarrow \mathfrak{c}\left( L\otimes V\right) \subseteq V\otimes L, \\
\mathfrak{cb}_{2} &=&\mathfrak{b}_{1}\mathfrak{c}_{2}\mathfrak{c}%
_{1}\Rightarrow \mathfrak{c}\left( V\otimes L\right) \subseteq L\otimes V.
\end{eqnarray*}%
Thus $L$ is a categorical subspace of $V.$ By Theorem \ref{teo: categorical}%
, we get that $L=0$ or $L=V$.
\end{proof}

\begin{proposition}
Let $V$ be an object in the monoidal category $_{H}^{H}\mathcal{YD}$ of
Yetter-Drinfeld modules over some Hopf algebra $H.$ Assume that $c_{V,V}$ is
a braiding of Hecke type of mark $\lambda $ and that $\lambda \neq 0,1.$
Then $V$ is simple in $_{H}^{H}\mathcal{YD}$.
\end{proposition}

\begin{proof}
Any subspace of $V$ in $_{H}^{H}\mathcal{YD}$ is categorical. We conclude by
Theorem \ref{teo: categorical}.
\end{proof}

\section{Trivial braided brackets}

Let $\mathfrak{b}$ be a braided bracket on $(V,\mathfrak{c})$. Our aim now
is to answer the following natural question: when is $\iota _{\mathfrak{c}{},%
\mathfrak{b}}:V\rightarrow U(V,\mathfrak{c},\mathfrak{b})$ injective?

\begin{proposition}
\label{pr: S izo grU}Let $\mathfrak{b}$ be a $\mathfrak{c}$-bracket on a
braided vector space $(V,\mathfrak{c})$ of Hecke type. If $\mathrm{gr}%
^{\prime }U(V,\mathfrak{c},\mathfrak{b})$ is the graded associated to the
standard filtration on $U(V,\mathfrak{c},\mathfrak{b}),$ then $\mathrm{gr}%
^{\prime }U(V,\mathfrak{c},\mathfrak{b})$ is a graded braided bialgebra and
there is a canonical morphism $\theta :S(V,\mathfrak{c})\rightarrow \mathrm{%
gr}^{\prime }U(V,\mathfrak{c},\mathfrak{b})$ of graded braided bialgebras.
Moreover $\theta $ is surjective.
\end{proposition}

\begin{proof}
Let $T^{\leq n}:=\bigoplus_{0\leq m\leq n}T^{m}$ and let $(U_{n}^{\prime
})_{n\in \mathbb{N}}$ be the standard filtration on $U:=U(V,\mathfrak{c},%
\mathfrak{b}).$ Let $\nabla _{U}$ and $\mathfrak{c}_{U}$ be the
multiplication and the braiding of $U:=U(V,\mathfrak{c},\mathfrak{b})$,
respectively$\mathfrak{.}$ If $T:=T(V,\mathfrak{c})$ and $\mathfrak{c}_{T}$
is the braiding of $T$ then the canonical projection $\pi _{U}:T\rightarrow
U $ is a morphism of braided bialgebras. Since $\mathfrak{c}_{T}(T^{\leq
n}\otimes T^{\leq m})\subseteq T^{\leq m}\otimes T^{\leq n}$ we deduce that $%
\mathfrak{c}_{U}(U_{n}^{\prime }\otimes U_{m}^{\prime })\subseteq
U_{m}^{\prime }\otimes U_{n}^{\prime },$ for any $n,m\in \mathbb{N}.$ Hence $%
\mathfrak{c}_{U}$ induces a braiding $\mathfrak{c}_{\mathrm{gr}^{\prime }U}:%
\mathrm{gr}^{\prime }U\otimes \mathrm{gr}^{\prime }U\rightarrow \mathrm{gr}%
^{\prime }U\otimes \mathrm{gr}^{\prime }U.$ The standard filtration is a
coalgebra filtration, as $\pi _{U}$ is a morphism of coalgebras, so $\mathrm{%
gr}^{\prime }U$ is a coalgebra. One can prove easily that, with respect to
this coalgebra structure, $\mathrm{gr}^{\prime }U$ becomes a graded braided
bialgebra.

Let $\lambda =\mathfrak{m}(\mathfrak{c}).$ We define $\theta
^{1}:V\rightarrow U_{1}^{\prime }/U_{0}^{\prime }$ by $\theta ^{1}=p\iota _{%
\mathfrak{c},\mathfrak{b}},$ where $p:U_{1}^{\prime }\rightarrow
U_{1}^{\prime }/U_{0}^{\prime }$ is the canonical projection. The image of $%
\theta ^{1}$ is included in the component of degree $1$ of $\mathrm{gr}%
^{\prime }U,$ so $\mathrm{Im\,}\theta ^{1}\subseteq P(\mathrm{gr}^{\prime
}U).$ Clearly $\theta ^{1}$ is a map of braided vector spaces. One can check
that $\nabla _{\mathrm{gr}^{\prime }U}\mathfrak{c}_{\mathrm{gr}^{\prime
}U}=\lambda \nabla _{\mathrm{gr}^{\prime }U}$ on $U_{1}^{\prime
}/U_{0}^{\prime }\otimes U_{1}^{\prime }/U_{0}^{\prime }$. By Proposition %
\ref{pr: universalitate S} (see also Remark \ref{re: universalitate S})
there is a unique morphism of graded braided bialgebras $\theta :S(V,%
\mathfrak{c})\rightarrow \mathrm{gr}^{\prime }U$ that lifts $\theta ^{1}$.
On the other hand, $\mathrm{gr}^{\prime }U(V,\mathfrak{c},\mathfrak{b})$ is
generated as an algebra by $U_{1}^{\prime }/U_{0}^{\prime }.$ Since $%
U_{1}^{\prime }/U_{0}^{\prime }$ is included into the image of $\theta ,$ we
conclude that $\theta $ is surjective.
\end{proof}

\begin{theorem}
\label{te: Lie algebras}Let $K$ be a field with $\mathrm{char}\,K\neq 2.$
Let $\mathfrak{b}$ be a $\mathfrak{c}$-bracket on a braided vector space $(V,%
\mathfrak{c})$ of Hecke-type of mark $\lambda \neq 0$ such that $%
(3)!_{\lambda }\neq 0$.

Assume that the $K$-linear map $\iota _{\mathfrak{c},\mathfrak{b}%
}:V\rightarrow U(V,\mathfrak{c},\mathfrak{b})$ is injective.

Then $\mathfrak{b}=0$ unless $\lambda =1.$ In this case $\mathfrak{b}$
fulfills%
\begin{equation}
\mathfrak{bc}=-\mathfrak{b}\text{\qquad and\qquad }\mathfrak{bb}_{1}(\mathrm{%
Id}_{V^{\otimes 3}}-\mathfrak{c}_{2}+\mathfrak{c}_{2}\mathfrak{c}_{1})=0
\label{form: Kharchenko}
\end{equation}%
i.e. $(V,\mathfrak{c},\mathfrak{b})$ is a generalized Lie algebra in the
sense of \cite{Gu}.
\end{theorem}

\begin{proof}
Denote by $T^{n}$ the $n$-th graded component of $T=T\left( V,c\right) $ and
set $T^{\leq n}:=\bigoplus_{0\leq m\leq n}T^{m}$ and $T^{\geq
n}:=\bigoplus_{m\geq n}T^{m}.$

Let $\gamma :=\mathfrak{c}-\lambda \mathrm{Id}_{T^{2}}-\mathfrak{b}$, let $%
F:=\mathrm{\mathrm{Im\,}}(\gamma )$ and $R:=\mathrm{\mathrm{Im\,}}\left(
\lambda \mathrm{Id}_{V^{\otimes 2}}-\mathfrak{c}\right) $. We denote the
component of degree $1$ of the map $\theta $ of Proposition \ref{pr: S izo
grU} by $\theta ^{1}:V\rightarrow U_{1}^{\prime }/U_{0}^{\prime }$.

Let $x\in \mathrm{Ker\,}\theta ^{1}.$ It follows that $\iota _{\mathfrak{c},%
\mathfrak{b}}(x)\in U_{0}^{\prime },$ so there is $a\in K$ such that $%
x-a1\in \left( F\right) .$ Since $(F)\subseteq T^{\geq 1}$ we get $a=0.$ It
results that $\iota _{\mathfrak{c},b}(x)=0,$ so $x=0.$ In conclusion, $%
\theta ^{1}$ is injective. But, in view of Theorem \ref{teo: S(V,c) type one}%
, $S(V,\mathfrak{c})$ is a strict coalgebra so that $P\left( S(V,\mathfrak{c}%
)\right) =S^{1}(V,\mathfrak{c})=V.$ Thus, since $\theta ^{1},$ the
restriction of $\theta $ to $V$, is injective, by \cite[Lemma 5.3.3]{Mo} it
follows that $\theta $ is injective too. Since $\theta $ is always
surjective, we conclude that $\theta $ is an isomorphism of braided
bialgebras.

The algebras $S(V,\mathfrak{c})$ and $U(V,\mathfrak{c},\mathfrak{b})$ are
the quotients of $T(V,\mathfrak{c})$ through the two-sided ideals generated
by $R$ and $F$, respectively. Set
\begin{equation*}
\zeta =(\lambda \mathrm{Id}_{V^{\otimes 3}}-\mathfrak{c}_{1})(\lambda ^{2}%
\mathrm{Id}_{V^{\otimes 3}}-\lambda \mathfrak{c}_{2}+\mathfrak{c}_{2}%
\mathfrak{c}_{1})\overset{\text{(\ref{ec: braid equation})}}{=}(\lambda
\mathrm{Id}_{V^{\otimes 3}}-\mathfrak{c}_{2})(\lambda ^{2}\mathrm{Id}%
_{V^{\otimes 3}}-\lambda \mathfrak{c}_{1}+\mathfrak{c}_{1}\mathfrak{c}_{2}).
\end{equation*}%
Since $(2)_{\lambda }\neq 0,$ one can easily see that $R=\{x\in T^{2}\mid
\mathfrak{c}(x)=-x\}$ and since $(3)!_{\lambda }\neq 0$ one gets that
\begin{equation}
(R\otimes V)\cap (V\otimes R)=\{x\otimes T^{3}(V)\mid \mathfrak{c}_{1}(x)=%
\mathfrak{c}_{2}(x)=-x\}=\mathrm{Im\,}\zeta .  \label{form: zeta}
\end{equation}%
Since the canonical map $\theta :S(V,\mathfrak{b})\rightarrow \mathrm{gr}%
^{\prime }U(V,\mathfrak{c},\mathfrak{b})$ is an isomorphism, by \cite[Lemma
0.4]{BG}, it follows that the following conditions are satisfied:%
\begin{gather}
F\cap T^{\leq 1}=0  \label{ec: injective 1} \\
\left( T^{\leq 1}\cdot F\cdot T^{\leq 1}\right) \cap T^{\leq 2}=F.
\label{ec: injective 2}
\end{gather}%
We claim (\ref{ec: injective 1}) implies
\begin{equation}
\mathfrak{bc}=-\mathfrak{b.}  \label{ec: antisymmetry}
\end{equation}%
In fact we have%
\begin{equation}
\gamma \left( \mathfrak{c}+\mathrm{Id}_{T^{2}}\right) =\left( \mathfrak{c}%
-\lambda \mathrm{Id}_{T^{2}}-\mathfrak{b}\right) \left( \mathfrak{c}+\mathrm{%
Id}_{T^{2}}\right) =-\mathfrak{b}\left( \mathfrak{c}+\mathrm{Id}%
_{T^{2}}\right)  \label{form: bella}
\end{equation}%
so that $\mathrm{Im}\left[ \mathfrak{b}\left( \mathfrak{c}+\mathrm{Id}%
_{T^{2}}\right) \right] =\mathrm{Im}\left[ \gamma \left( \mathfrak{c}+%
\mathrm{Id}_{T^{2}}\right) \right] .$ We will prove that $\mathrm{Im}\left[
\gamma \left( \mathfrak{c}+\mathrm{Id}_{T^{2}}\right) \right] =F\cap T^{\leq
1},$ from which the conclusion will follow.

$\subseteq )$ It follows by (\ref{form: bella}).

$\supseteq )$ Let $y\in F\cap T^{\leq 1}.$ Then there is $x$ such that $%
y=\gamma \left( x\right) =\mathfrak{c}\left( x\right) -\lambda x-\mathfrak{b}%
\left( x\right) .$ Since $y,\mathfrak{b}\left( x\right) \in T^{\leq 1}$ and $%
\mathfrak{c}\left( x\right) ,x\in V\otimes V,$ it results $\mathfrak{c}%
(x)=\lambda x$ and $y=-\mathfrak{b}(x).$ Thus

$\gamma \left( \mathfrak{c}+\mathrm{Id}_{T^{2}}\right) \left( x\right)
=\left( 2\right) _{\lambda }\gamma \left( x\right) =\left( 2\right)
_{\lambda }y.$ Since $\left( 2\right) _{\lambda }\neq 0,$ we get $y\in
\mathrm{Im}\left[ \gamma \left( \mathfrak{c}+\mathrm{Id}_{T^{2}}\right) %
\right] .$

Now, we define $\alpha :=-\left( 2\right) _{\lambda }^{-1}\mathfrak{b}|_{R}.$
From (\ref{ec: antisymmetry}), we deduce that
\begin{equation*}
\left( \mathrm{Id}_{T^{2}}-\alpha \right) \left( \mathfrak{c}-\lambda
\mathrm{Id}_{T^{2}}\right) =\mathfrak{c}-\lambda \mathrm{Id}_{T^{2}}+\frac{1%
}{\left( 2\right) _{\lambda }}\mathfrak{b}\left( \mathfrak{c}-\lambda
\mathrm{Id}_{T^{2}}\right) =\gamma
\end{equation*}%
so that $F=\mathrm{Im}\gamma =\mathrm{Im}\left[ \left( \mathrm{Id}%
_{T^{2}}-\alpha \right) \left( \mathfrak{c}-\lambda \mathrm{Id}%
_{T^{2}}\right) \right] =\{x-\alpha (x)\mid x\in R\}.$ By \cite[Lemma 3.3]%
{BG} (\ref{ec: injective 2}) implies that $\alpha $ satisfies the following
two conditions:
\begin{eqnarray*}
(\alpha \otimes V)(x)-(V\otimes \alpha )(x) &\in &R,\qquad \forall x\in
(R\otimes V)\cap (V\otimes R), \\
\alpha (\alpha \otimes V-V\otimes \alpha )(x) &=&0,\qquad \forall x\in
(R\otimes V)\cap (V\otimes R).
\end{eqnarray*}%
The second property is equivalent to the fact that $\mathfrak{b}(\mathfrak{b}%
_{1}-\mathfrak{b}_{2})=0$ on $\mathrm{Im\,}\zeta ,$ which at its turn is
equivalent to
\begin{equation}
\mathfrak{bb}_{1}\zeta =\mathfrak{bb}_{2}\zeta .  \label{ec: Jacoby}
\end{equation}%
Let us prove that $\mathfrak{b}$ satisfies the following conditions
\begin{equation}
\mathfrak{bb}_{1}(\lambda ^{2}\mathrm{Id}_{V^{\otimes 3}}-\lambda \mathfrak{c%
}_{2}+\mathfrak{c}_{2}\mathfrak{c}_{1})=0,  \label{ec: Jacoby 1}
\end{equation}%
\begin{equation}
\mathfrak{bb}_{2}(\lambda ^{2}\mathrm{Id}_{V^{\otimes 3}}-\lambda \mathfrak{c%
}_{1}+\mathfrak{c}_{1}\mathfrak{c}_{2})=0.  \label{ec: Jacoby 2}
\end{equation}%
We have
\begin{equation*}
\mathfrak{bb}_{2}\zeta \overset{\text{(\ref{form: zeta})}}{=}\mathfrak{bb}%
_{2}\mathfrak{c}_{1}\mathfrak{c}_{2}\zeta \overset{\text{(\ref{ec: b
compatibilitate})}}{=}\mathfrak{b\mathfrak{c}b}_{1}\zeta \overset{\text{(\ref%
{ec: antisymmetry})}}{=}-\mathfrak{bb}_{1}\zeta
\end{equation*}%
In view of (\ref{ec: Jacoby}), we obtain $\mathfrak{bb}_{2}\zeta =0=%
\mathfrak{bb}_{1}\zeta $ as $\mathrm{char}\left( K\right) \neq 2.$ We have%
\begin{equation*}
0=\mathfrak{bb}_{1}\zeta =\mathfrak{bb}_{1}(\lambda \mathrm{Id}_{V^{\otimes
3}}-\mathfrak{c}_{1})(\lambda ^{2}\mathrm{Id}_{V^{\otimes 3}}-\lambda
\mathfrak{c}_{2}+\mathfrak{c}_{2}\mathfrak{c}_{1})\overset{\text{(\ref{ec:
antisymmetry})}}{=}(2)_{\lambda }\mathfrak{bb}_{1}(\lambda ^{2}\mathrm{Id}%
_{V^{\otimes 3}}-\lambda \mathfrak{c}_{2}+\mathfrak{c}_{2}\mathfrak{c}_{1}).
\end{equation*}%
and%
\begin{equation*}
0=\mathfrak{bb}_{2}\zeta =\mathfrak{bb}_{2}(\lambda \mathrm{Id}_{V^{\otimes
3}}-\mathfrak{c}_{2})(\lambda ^{2}\mathrm{Id}_{V^{\otimes 3}}-\lambda
\mathfrak{c}_{1}+\mathfrak{c}_{1}\mathfrak{c}_{2})\overset{\text{(\ref{ec:
antisymmetry})}}{=}(2)_{\lambda }\mathfrak{bb}_{2}(\lambda ^{2}\mathrm{Id}%
_{V^{\otimes 3}}-\lambda \mathfrak{c}_{1}+\mathfrak{c}_{1}\mathfrak{c}_{2}).
\end{equation*}%
so that $\mathfrak{b}$ satisfies (\ref{ec: Jacoby 1}) and (\ref{ec: Jacoby 2}%
). We have%
\begin{equation*}
0\overset{\text{(\ref{ec: b compatibilitate})}}{=}\mathfrak{b}\left[
\mathfrak{b}_{1}\mathfrak{c}_{2}\mathfrak{c}_{1}-\mathfrak{cb}_{2}\right] =%
\mathfrak{bb}_{1}\mathfrak{c}_{2}\mathfrak{c}_{1}-\mathfrak{bcb}_{2}\overset{%
\text{(\ref{ec: antisymmetry})}}{=}\mathfrak{bb}_{1}\mathfrak{c}_{2}%
\mathfrak{c}_{1}+\mathfrak{bb}_{2}\overset{\text{(\ref{ec: Jacoby 1})}}{=}%
\mathfrak{bb}_{1}\left( -\lambda ^{2}\mathrm{Id}_{V^{\otimes 3}}+\lambda
\mathfrak{c}_{2}\right) +\mathfrak{bb}_{2}
\end{equation*}%
so that $\mathfrak{bb}_{2}=\lambda \mathfrak{bb}_{1}\left( \lambda \mathrm{Id%
}_{V^{\otimes 3}}-\mathfrak{c}_{2}\right) .$ Similarly using $\mathfrak{cb}%
_{1}=\mathfrak{b}_{2}\mathfrak{c}_{1}\mathfrak{c}_{2}$ in (\ref{ec: b
compatibilitate}), (\ref{ec: antisymmetry}) and (\ref{ec: Jacoby 2}) we get $%
\mathfrak{bb}_{1}=\lambda \mathfrak{bb}_{2}\left( \lambda \mathrm{Id}%
_{V^{\otimes 3}}-\mathfrak{c}_{1}\right) .$ Using these formulas we obtain%
\begin{eqnarray*}
\mathfrak{bb}_{2} &=&\lambda \mathfrak{bb}_{1}\left( \lambda \mathrm{Id}%
_{V^{\otimes 3}}-\mathfrak{c}_{2}\right) =\lambda \lambda \mathfrak{bb}%
_{2}\left( \lambda \mathrm{Id}_{V^{\otimes 3}}-\mathfrak{c}_{1}\right)
\left( \lambda \mathrm{Id}_{V^{\otimes 3}}-\mathfrak{c}_{2}\right) \\
&=&\lambda ^{2}\mathfrak{bb}_{2}\left( \lambda ^{2}\mathrm{Id}_{V^{\otimes
3}}-\lambda \mathfrak{c}_{1}+\mathfrak{c}_{1}\mathfrak{c}_{2}-\lambda
\mathfrak{c}_{2}\right) \overset{\text{(\ref{ec: Jacoby 2})}}{=}-\lambda ^{3}%
\mathfrak{bb}_{2}\mathfrak{c}_{2}\overset{\text{(\ref{ec: antisymmetry})}}{=}%
\lambda ^{3}\mathfrak{bb}_{2}
\end{eqnarray*}%
so that $\left( \lambda ^{3}-1\right) \mathfrak{bb}_{2}=0.$ Assume $%
\mathfrak{b}$ is not zero. By Proposition \ref{pro: triviality}, if $\lambda
\neq 0,1$, we get that $\mathfrak{b}$ is surjective. Thus $\mathfrak{bb}_{2}=%
\mathfrak{b}\left( V\otimes \mathfrak{b}\right) $ is surjective too and
hence $\lambda ^{3}=1,$ contradicting $(3)_{\lambda }\neq 0.$
\end{proof}

\begin{remarks}
Concerning the converse of Theorem \ref{te: Lie algebras}, let us
note that
if $\mathfrak{b}=0$ then the map $\iota _{\mathfrak{c},\mathfrak{b}%
}:V\rightarrow U(V,\mathfrak{c},\mathfrak{b})=S(V,\mathfrak{c})$ is clearly
injective. On the other hand, let $K$ be a field of characteristic zero.
Given a $\mathfrak{c}$-bracket $\mathfrak{b}$ on a braided vector space $(V,%
\mathfrak{c})$ where $\mathfrak{c}^{2}=\mathrm{Id}_{V\otimes V}$ (i.e. $(V,%
\mathfrak{c})$ of Hecke-type of regular mark $1$), then the canonical map $%
\iota _{\mathfrak{c},\mathfrak{b}}:V\rightarrow U(V,\mathfrak{c},\mathfrak{b}%
)$ is injective whenever $\mathfrak{b}$ fulfills (\ref{form: Kharchenko})
(see \cite[Theorem 5.2]{Kharchenko- connected}).
\end{remarks}

\begin{example}
\cite{Masuoka} Let $K$ be a field with $\mathrm{char}\left( K\right) =2$.
Let $V=Kx$ and let $\mathfrak{c}:V\otimes V\rightarrow V\otimes V,\mathfrak{c%
}=\mathrm{Id}_{V\otimes V}.$ Define $\mathfrak{b}:V\otimes V\rightarrow V$
by $\mathfrak{b}\left( x\otimes x\right) =ax$ for some $a\in K\backslash
\left\{ 0\right\} .$

Assume there exists $\lambda \in K\backslash \left\{ 0,1\right\} $ such that
$(3)!_{\lambda }\neq 0$ i.e. such that $\lambda $ is not a primitive third
root of unity. Clearly $\left( \mathfrak{c}+\mathrm{Id}_{V\otimes V}\right)
\left( \mathfrak{c}-\lambda \mathrm{Id}_{V\otimes V}\right) =0$ so that $%
\mathfrak{c}$ is of Hecke-type of mark $\lambda $. Moreover $\mathfrak{b}$
is a $\mathfrak{c}$-bracket on the braided vector space $\left( V,\mathfrak{c%
}\right) $. Thus we can consider the universal enveloping algebra%
\begin{equation*}
U=U(V,\mathfrak{c},\mathfrak{b}):=\frac{T(V,{\mathfrak{c}})}{\left(
\mathfrak{c}(z)-\lambda z-\mathfrak{b}(z)\mid z\in V\otimes V\right) }\simeq
\frac{K\left[ X\right] }{\left( \left( 1-\lambda \right) X^{2}-aX\right) }.
\end{equation*}%
The canonical map $\iota _{\mathfrak{c},\mathfrak{b}}:V\rightarrow U(V,%
\mathfrak{c},\mathfrak{b})$ is clearly injective. Nevertheless $\mathfrak{b}%
\neq 0.$
\end{example}

\section{A Milnor-Moore type theorem for braided bialgebras}

In this section we prove the main result of this paper, Theorem \ref{te: MM}%
, which represents a variant of Milnor-Moore Theorem for braided bialgebras.
Then we deduce some consequences of this theorem, including applications to
certain classes of bialgebras in braided categories.

\begin{definition}
Let $\left( A,\mathfrak{c}_{A}\right) $ be a connected braided bialgebra.
Let $P:=P(A)$. The braiding $\mathfrak{c}_{P}=\mathfrak{c}_{A}|_{P\otimes P}$
will be called the\textbf{\ infinitesimal braiding }of $A$.
\end{definition}

\begin{remarks}
\label{rem: infinitesimal}Let $\left( A,\mathfrak{c}_{A}\right) $ be a
connected braided bialgebra. Let $P:=P(A)$ and let $\mathfrak{c}_{P}\ $be
the infinitesimal braiding of $A.$ If $\mathrm{gr}\,A$ denotes the graded
associated with respect to the coradical filtration, then $\mathrm{gr}\,A$
is strictly graded. Thus
\begin{equation*}
P(\mathrm{gr}\,A)=\mathrm{gr}^{1}A\simeq P(A)=P.
\end{equation*}%
Through this identification, $\mathfrak{c}_{\mathrm{gr}\,A}^{1,1}$ is equal
to $\mathfrak{c}_{P}$. In conclusion the infinitesimal braiding of $\mathrm{%
gr}\,A$ is the infinitesimal braiding $\mathfrak{c}_{P}$ of $A$.
\end{remarks}

\begin{definition}
Let $\left( A,\mathfrak{c}_{A}\right) $ be a connected braided bialgebra and
let $P:=P(A)$. The component $\Delta _{\mathrm{gr}\,A}^{1,1}:A_{2}/A_{1}%
\rightarrow A_{1}/A_{0}\otimes A_{1}/A_{0}=P\otimes P$ is called the \emph{%
infinitesimal comultiplication of} $A$.

Let $\mathfrak{c}_{P}=\mathfrak{c}_{A}|_{P\otimes P}$ and let $\lambda \in
K^{\ast }$. We will say that $\Delta _{\mathrm{gr}\,A}^{1,1}$ is $\lambda $-%
\emph{cocommutative} if $\mathfrak{c}_{P}\circ \Delta _{\mathrm{gr}%
\,A}^{1,1}=\lambda \Delta _{\mathrm{gr}\,A}^{1,1}$, that is we have:
\begin{equation}
\mathfrak{c}_{\mathrm{gr}\,A}^{1,1}\circ \Delta _{\mathrm{gr}%
\,A}^{1,1}=\lambda \Delta _{\mathrm{gr}\,A}^{1,1}.  \label{ec: Delta}
\end{equation}
\end{definition}

\begin{proposition}
\label{pr: universalitate U2}Let $K$ be a field with $\mathrm{char}\,K\neq
2. $ Let $A$ be a connected braided bialgebra and assume that its
infinitesimal braiding is of Hecke-type of mark $\lambda \neq 0,1$ such that
$(3)!_{\lambda }\neq 0$. Let $P$ be the space of primitive elements of $A$
and let $\mathfrak{b}_{P}=\nabla (\mathfrak{c}_{P}-\lambda \mathrm{Id}%
_{P^{\otimes 2}})|_{P\otimes P}$ be the $\mathfrak{c}_{P}$-bracket on the
braided vector space $(P,\mathfrak{c}_{P})$ defined in Proposition \ref{pr:
universalitate U}. Then $\mathfrak{b}_{P}=0$.

Let $f:(V,\mathfrak{c},\mathfrak{b})\rightarrow (P,\mathfrak{c}_{P},0)$ be a
morphism of braided brackets and assume that $\mathfrak{c}$ is a braiding of
Hecke-type with mark $\lambda .$ Then there is a unique morphism of braided
bialgebras $\widetilde{f}:U(V,\mathfrak{c},\mathfrak{b})\rightarrow A$ that
lifts $f$.
\end{proposition}

\begin{proof}
By Proposition \ref{pr: universalitate U}(b) it follows that $\mathfrak{b}%
_{P}$ is a $\mathfrak{c}_{P}$-bracket on $(P,\mathfrak{c}_{P})$, hence we
can apply the universal property of $U:=U(P,\mathfrak{c}_{P},\mathfrak{b}%
_{P})$. There is a unique morphism of braided bialgebras $\phi
_{A}:U\rightarrow A$ that lifts $\mathrm{Id}_{P}$. Observe that the
canonical map $\iota _{\mathfrak{c},\mathfrak{b}}:P\rightarrow U$ is
injective, as $\phi _{A}\iota _{\mathfrak{c},\mathfrak{b}}$ is the inclusion
of $P$ into $A$. Now apply Theorem \ref{te: Lie algebras} to obtain that $%
\mathfrak{b}_{P}=0$. The last part follows by \ref{pr: universalitate U}.
\end{proof}

\begin{theorem}
\label{te: MM}Let $K$ be a field with $\mathrm{char}\,K\neq 2.$ Let $A$ be a
connected braided bialgebra and assume that its infinitesimal braiding $%
\mathfrak{c}_{P\left( A\right) }$ is $\lambda $-cocommutative for some
regular element $\lambda \neq 0$ in $K$. Then

\begin{itemize}
\item the infinitesimal braiding of $A$ is of Hecke-type of mark $\lambda $
and

\item $A$ is isomorphic as a braided bialgebra to the symmetric algebra $%
S\left( P\left( A\right) ,\mathfrak{c}_{P\left( A\right) }\right) $ of $%
\left( P\left( A\right) ,\mathfrak{c}_{P\left( A\right) }\right) $ whenever $%
\lambda \neq 1$.
\end{itemize}
\end{theorem}

\begin{proof}
Let $B:=\mathrm{gr}A$. Clearly $B$ is strongly $%
%TCIMACRO{\U{2115} }%
%BeginExpansion
\mathbb{N}
%EndExpansion
$-graded as a coalgebra. By assumption the infinitesimal braiding of $A$ is $%
\lambda $-cocommutative and hence the same holds for $B$ i.e. $\left(
c_{B}^{1,1}-\lambda \mathrm{Id}_{B^{2}}\right) \Delta _{B}^{1,1}=0$. Since $B
$ is also $0$-connected, by Theorem \ref{teo: aureo}, $B$ is a bialgebra of
type one and $\mathfrak{c}_{B}^{1,1}$ is a braiding of Hecke-type of mark $%
\lambda $. In particular the infinitesimal braiding of $A$ is of Hecke-type
of mark $\lambda $ and $B$ is generated as an algebra by $B^{1}\ $so that $A$
is generated as a $K$-algebra by $P=P\left( A\right) =B^{1}.$ Therefore, the
canonical braided bialgebra homomorphism $f:U\left( P,c_{P},\mathfrak{b}%
_{P}\right) \rightarrow A$, arising by the universal property of the
universal enveloping algebra, is surjective. Assume $\lambda \neq 1$. By
Proposition \ref{pr: universalitate U2}, $\mathfrak{b}_{P}=0$ hence $U\left(
P,c_{P},\mathfrak{b}_{P}\right) =S\left( P,c_{P}\right) .$ On the other
hand, by Theorem \ref{teo: S(V,c) type one}, $P$ is the primitive part of $%
S\left( P,c_{P}\right) $ and the restriction of $f$ to $P$ is injective so
that $f$ is injective by \cite[Lemma 5.3.3]{Mo}. In conclusion $f$ is an
isomorphism.
\end{proof}

\begin{remark}
\label{rem: MM}With hypothesis of Theorem \ref{te: MM}, if $\lambda =1$ then
$A$ is isomorphic as a braided bialgebra to the universal enveloping algebra
$U\left( P\left( A\right) ,\mathfrak{c}_{P\left( A\right) },\mathfrak{b}%
_{P\left( A\right) }\right) $ of $\left( P\left( A\right) ,\mathfrak{c}%
_{P\left( A\right) },\mathfrak{b}_{P\left( A\right) }\right) .$ In fact
regularity of $\lambda $ in this case means $\mathrm{char}\left( K\right)
=0. $ By Theorem \ref{te: Lie algebras} $\mathfrak{bc}=-\mathfrak{b}$ and $%
\mathfrak{bb}_{1}(\mathrm{Id}_{V^{\otimes 3}}-\mathfrak{c}_{2}+\mathfrak{c}%
_{2}\mathfrak{c}_{1})=0$ so that \cite[Theorem 6.1]{Kharchenko- connected}
applies.
\end{remark}

\begin{corollary}
\label{co:graded S(V)}Let $K$ be a field with $\mathrm{char}\,K\neq 2.$ Let $%
(V,\mathfrak{c})$ be a braided vector space such that $\mathfrak{c}$ is a
braiding of Hecke-type of regular mark $\lambda \neq 0,1$. Let $A$ be a
braided bialgebra such that ${\mathrm{gr}}\,A$ is isomorphic as a braided
bialgebra to $S(V,\mathfrak{c})$ then $A$ is isomorphic to the symmetric
algebra $S(V,\mathfrak{c})$ of $(V,\mathfrak{c}).$
\end{corollary}

\begin{proof}
Obviously the infinitesimal comultiplication of $S(V,\mathfrak{c})$ is $%
\lambda $-cocommutative and, by Proposition \ref{pr: Sconn}, $S(V,\mathfrak{c%
})$ is connected. Thus $\mathrm{gr}\,A$ has the same properties. Hence $A$
itself is connected and by Remark \ref{rem: infinitesimal}, the
infinitesimal braiding of $A$ is $\lambda $-cocommutative. We conclude by
applying Theorem \ref{te: MM}.
\end{proof}

\noindent%
\begin{minipage}[t]{15cm}\vspace*{2mm}\sc \footnotesize
University of Ferrara, Department of Mathematics, Via Machiavelli
35, I-44100, Ferrara,  Italy. \\
{\small\it email:} {\small\ttfamily \rm alessandro.ardizzoni@unife.it}\medskip\\%
University of Ferrara, Department of Mathematics, Via Machiavelli
35, I-44100, Ferrara, Italy\\
{\small\it email:} {\small\ttfamily \rm men@unife.it}\medskip\\%
University of Bucharest, Faculty of Mathematics, Academiei 14,
RO-010014, Bucharest, Romania. \\
{\small\it email:} {\small\ttfamily \rm dstefan@al.math.unibuc.ro}
\end{minipage}

\end{document}